\documentclass[12pt, reqno]{amsart}
\usepackage{amssymb, amsthm, amsmath, amsfonts, amsaddr}
\usepackage{array, epsfig}
\usepackage{bbm}
\usepackage{pgf}

  \evensidemargin
\oddsidemargin
\begin{document}

\newcommand{\be}{\begin{equation}}
\newcommand{\ee}{\end{equation}}
\newcommand{\bea}{\begin{eqnarray}}
\newcommand{\eea}{\end{eqnarray}}
\newcommand{\beaa}{\begin{eqnarray*}}
\newcommand{\eeaa}{\end{eqnarray*}}

\renewcommand{\proofname}{\bf Proof}
\newtheorem*{rem*}{Remark}
\newtheorem*{cor*}{Corollary}
\newtheorem{cor}{Corollary}
\newtheorem{prop}{Proposition}
\newtheorem{lem}{Lemma}
\newtheorem*{lem1'}{Lemma $\mathbf{1^\prime}$}
\newtheorem{theo}{Theorem}
\newfont{\zapf}{pzcmi}

\def\R{\mathbb{R}}
\def\Z{\mathbb{Z}}
\def\N{\mathbb{N}}
\def\E{\mathbb{E}}
\def\P{\mathbb{P}}
\def\V{\mathbb{D}}
\def\I{\mathbbm{1}}
\newcommand{\D}{\hbox{\zapf D}}
\newcommand{\Pt}{\widetilde{\mathbb{P}}}
\newcommand{\Et}{\widetilde{\mathbb{E}}}

\newcommand{\bt}{\begin{theo}}
\newcommand{\et}{\end{theo}}
\newcommand{\bl}{\begin{lem}}
\newcommand{\el}{\end{lem}}
\newcommand{\bc}{\begin{cor*}}
\newcommand{\ec}{\end{cor*}}
\newcommand{\br}{\begin{rem*}}
\newcommand{\er}{\end{rem*}}
\newcommand{\bp}{\begin{proof}}
\newcommand{\ep}{\end{proof}}
\newcommand{\bes}{\begin{ex}}
\newcommand{\ees}{\end{ex}}

\title[Limit theorems for random walks that avoid bounded sets]{Limit theorems for random walks that avoid bounded sets, with applications to the largest gap problem}
\author{Vladislav Vysotsky}
\address{Arizona State University, St.Petersburg Department of Steklov Mathematical Institute, Chebyshev Laboratory at St.Petersburg State University}
\email{vysotsky@asu.edu, vysotsky@pdmi.ras.ru}
\thanks{This research is supported in part by the Chebyshev Laboratory  (Department of Mathematics and Mechanics, St. Petersburg State University) under RF Government Grant 11.G34.31.0026, JSC ``Gazprom Neft'', Grants 13-01-00256 and 14-01-31375 by RFBR, and NSh-1216.2012.1}

\begin{abstract}
Consider a centred random walk in dimension one with a positive finite variance $\sigma^2$, and let $\tau_B$ be the hitting time for a bounded Borel set $B$ with a non-empty interior. We prove the asymptotic $\P_x(\tau_B > n) \sim \sqrt{2 / \pi} \sigma^{-1} V_B(x) n^{-1/2}$ and provide an explicit formula for the limit $V_B$ as a function of the initial position $x$ of the walk. We also give a functional limit theorem for the walk conditioned to avoid $B$ by the time $n$. As a main application, consider the case that $B$ is an interval and study the size of the largest gap $G_n$ (maximal spacing) within the range of the walk by the time $n$. We prove a limit theorem for $G_n$, which is shown to be of the constant order, and describe its limit distribution. In addition, we prove an analogous result for the number of non-visited sites within the range of an integer-valued random walk.


\end{abstract}

\maketitle

\section{Introduction and results}
\subsection{Introduction}

Let $S_n = x + X_1+ \dots + X_n$ be a centred random walk in dimension one with a positive finite variance $\sigma^2$. Denote $\tau_B := \min\{ k \ge 0: S_k \in B\}$ the hitting time of a Borel set $B$. The goal of this paper is to find the asymptotic of $\P_x(\tau_B >n)$ as $n \to \infty$ in the case that $B$ is bounded. Here $\P_x$ denotes the distribution of the walk starting at $S_0= x$.

Recall that the distribution of $X_1$ is called $\lambda$-arithmetic if $\lambda >0$ is the maximal number such that $\P(X_1 \in \lambda \Z) =1$, and is called non-arithmetic if such a number does not exist. We set the state space $M$ of the walk to be $\lambda \Z$ for $\lambda$-arithmetic
walks and $\R$ for non-arithmetic walks. We will consider the starting points $x \in M$ only; we stress that this does not result in any loss of generality. To avoid trivialities, we {\it always} assume that $Int_M(B)$ is non-empty, where in the arithmetic case the interior is taken under the discrete topology.


Kesten and Spitzer \cite{KestenSpitzer} (Theorem 4a)  proved that for {\it any} 1-arithmetic random walk, for any finite non-empty $B \subset \Z$ and $x \in \Z$ there exists the finite limit of the ratios
\be \label{ratio limit}
\lim_{n \to \infty} \frac{\P_x(\tilde{\tau}_B > n)}{\P_0(\tilde{\tau}_{\{0\}} > n)} =: g_B(x),
\ee
where $\tilde{\tau}_B := \min\{ k \ge 1: S_k \in B\}$ (of course $\tilde{\tau}_B=\tau_B$ for $x \notin B$). This is true no matter if the walk is transient or recurrent, and also holds in dimension two. Moreover, \cite{KestenSpitzer} gives a formula for the limit function, which under our assumptions reads as
\be \label{KS g=}
g_B(x) =  \lim_{y \to \infty} \frac12 \E_x  \sum_{k=0}^{\tilde{\tau}_B} \I_{\{|X_k| = y \}}.\\
\ee

The proof of \eqref{ratio limit} in~\cite{KestenSpitzer} is by induction in the number of elements of $B$. Essentially, the problem simplifies to the study of $\P_x(\tilde{\tau}_{\{0\}} > n)$, which can be conducted using the theory of recurrent events. In addition, the theory easily implies the asymptotic of $\P_0(\tilde{\tau}_{\{0\}} > n)$. Under our assumptions this gives the asymptotic
\be \label{KS}
\lim_{n \to \infty} \sqrt{n} \P_x(\tau_B > n) = \sqrt{\frac{2}{\pi}} \sigma a_B(x)
\ee
with $a_B(x) := g_B(x) \I_{B^c}(x)$, $x \in \Z$.

Notice that the limit function $a_B(x)$ is harmonic for the random  walk killed as it enters $B$, in the sense that $\E_x a_B(\bar{S}_1) = a_B(x)$ for any $x$, where $\bar{S}_n:=S_{n \wedge \tau_B}$ denotes the killed walk, which is a Markov chain. Since in this paper we do not need to consider $x \in B$, let us rewrite the above as $\E_x a_B(S_1) = a_B(x)$ for $x \notin B$ and $a_B(x)=0$ on $B$. The latter ``boundary condition'' is very unnatural for harmonic analysis, and for example, no uniqueness theorem is proved for non-negative harmonic functions of this type.

Of course the approach of Kesten and Spitzer is helpless in the general case of non-arithmetic walks, which is of the main interest for this paper. We will use an entirely different approach that forces us to work in dimension one and assume the finiteness of variance. Although we believe that the proposed methodology provides a right tool to study the problem for random walks with infinite variance, this case is totally different from the one presented here, and there is absolutely no universality as in the proof by Kesten and Spitzer.


Our main result extends \eqref{KS} from arithmetic to general random walks. Also we provide an explicit expression for the limit function as opposed to the one in \eqref{KS g=}, which essentially replaces one limit by other and from this perspective is very implicit. With no surprise, the order of the tail asymptotic is the same.

Let us introduce the notation needed to define this limit. We will relate it to the well-studied tail behaviour of the hitting times for half-lines. We start considering the basic case that $B = (-d, d)$ is an interval. Here we think of $d > 0$ as of a fixed parameter and introduce the special notation $$p_n^{(d)}(x):= p_n(x):=\P_x(\tau_{(-d,d)} > n).$$ Define the moments $T_k$ when the random walk jumps from $B^c$ to $B$ as follows: $T_0:=0$ and
$$T_{k+1} := \min\{ n> T_k: S_{n-1} \ge d, S_n < d \} \, \wedge \, \min\{ n> T_k: S_{n-1} \le -d, S_n > -d\}$$
for $k \ge 0$, and denote $H_k:=S_{T_k}$; recall that here we consider only $|x| \ge d$. Now $\tau_{(-d,d)} = T_\kappa$, where $\kappa:= \min(k \ge 1: |H_k| < d)$.  In some cases it is easier to work under one probability $\P$ rather than $\P_x$ with varying $x$. For this purpose we consider a walk that starts at zero and avoids the set $B-x$, while the variables introduced above are redefined accordingly and denoted as $T_1(x), H_1(x),$ etc. All the given definitions clearly depend on the specific choice of $B$ (i.e. $d$ in the current case). Finally, put $H_\pm:=S_{\tau_{\pm(0, \infty)}}$.

If $x \ge d$, then $T_1$ clearly is the moment of the first entrance to the half-line $(-\infty, d)$. Similarly, if $x \le -d$, then $T_1$ is the time of the first entrance to $(-d, \infty)$. The tail of the hitting time for a half-line is well understood and its asymptotic can be written in the following form that conveniently puts together two options of starting to the right and to the left of $B= (-d,d)$: it holds that
\be \label{half-tail}
\P_x( T_1 >n) \sim \frac{\sqrt{2} |x - \E_x S_{T_1}|}{\sigma \sqrt{\pi n} }, \quad n \to \infty, \qquad x \in M \setminus Conv(B).
\ee

Denote $U_d(x):=|x - \E_x S_{T_1}|$ for $|x| \ge d$, $U_+(y):= y - \E_y H_-$ for $y \ge 0$, $U_-(y):= \E_y H_+ - y$ for $y \le 0$. Of course it holds that $U_d(x) = U_{sgn(x)}(x - sgn(x) d)$ for any $|x| \ge d$. Since the walk $S_n$ is centred, the function $U_+(y)$ with $y \ge 0$ is harmonic for the walk $S_n$ killed as it enters $(-\infty, 0)$, in the sense that it solves the equation $U_+(y) = \E_y U(S_1) \I_{\{\tau_{(-\infty, 0)} > 1\}}$ for $y \ge 0$. The analogous is true for $U_-(y)$. Since $\E |H_\pm| < \infty$, the Wald identity implies that $U_+(x)$ and $U_-(x)$ are proportional to the renewal functions for the renewal processes generated by the strictly ascending and descending ladder heights of $S_n$, respectively.

{\it To our surprise and to the best of our knowledge, there is no appropriate reference for \eqref{half-tail}}. Moreover, it was never stated in such a comprehensive form. This is why we discuss \eqref{half-tail} here in detail. Feller~\cite{Feller} (Sec. XII.7, 8 and XVIII.5) is a classical source, which considers positive and negative excursions that start at zero ($x = d$ and $x=-d$, respectively, in our setting). Note that from the standard approach of Sec. XII in~\cite{Feller} it is not clear how to combine the tail asymptotics for positive and negative excursions into a single equation, while \eqref{half-tail}  easily follows from the expressions for the constants (see Sec. XVIII) given in terms of $\E_0 S_{\tau_\pm}$. Further, Bingham et al.~\cite{BGT} (Theorem 8.9.12) states that \eqref{half-tail} is true for any {\it continuity point} $x$ of $U_d(x)$. This result and Lemma 2 by Bertoin and Doney~\cite{BertoinDoney} combined with the right-continuity of $U_d(x)$ and monotonicity in $x$ of $\P_x( T_1 >n)$ easily imply \eqref{half-tail} for {\it any} $x$. Moreover, we will discuss a uniform version of \eqref{half-tail} for $x=o(\sqrt{n})$, and provide an analogous result for the hitting times of closed half-lines, see Lemmae~\ref{LEM UNIFORM} and \ref{LEM UNIFORM'} in Sec.~\ref{SEC PROOFS MAIN} below.

\subsection{Main results}
Define $$V_d(x):= \E_x \sum_{i=1}^{\kappa} \bigl| S_{T_i} - S_{T_{i-1}} \bigr| = \E_x \sum_{i=1}^{\kappa} \bigl| H_{i} - H_{i-1} \bigr|,$$ where of course $V_d(x)\equiv 0$ on $(-d,d)$. We are ready to state the first result of the paper.

\bt
\label{THM EQUIV}
Let $S_n$ be a centred random walk with positive finite variance $\sigma^2$. Then for any $d >0$ and any $x \in M$ we have that
\be \label{main}
p_n(x) \sim \frac{\sqrt{2} V_d(x)}{\sigma \sqrt{\pi n}}, \quad n \to \infty.
\ee
The function $V_d(x)$ is harmonic for the walk $S_n$ killed as it enters $(-d,d)$, that is $V_d(x) = \E_x V_d(S_1)$ for $|x| \ge d$. It satisfies $0 < U_d(x) \le V_d(x) < \infty$ for $|x| \ge d$, $V_d(x) \sim |x|$ as $|x| \to \infty$, and $V_d(\pm(d + y)) - U_d(\pm (d +y)) \to 0$ as $d \to \infty$ for any fixed $y \ge 0$.

Moreover, \eqref{main} holds uniformly for $x =o(\sqrt{n})$, that is over $x \in [-a_n \sqrt{n}, a_n \sqrt{n}] \cap M$ for any fixed sequence $a_n \to 0$ such that $a_n \sqrt{n} \to \infty$.
\et

For an heuristic explanation, it is instructive to think that the walk that avoids an interval may perform few jumps over the interval at the small values of time and then typically stays to one side of $(-d,d)$ up to time $n$.

Now let $B$ is a general Borel bounded set. Denote $$r := \sup\{x: x \in B \}, \quad l:= \inf \{x: x \in B \}, \quad B_-:=(-\infty, l] \setminus B, \quad B_+:= [r, \infty) \setminus B.$$ We modify the definition of $T_k$ as follows since here we need to consider {\it all} the moments when the walk jumps over the edges $l$ and $r$ of $B$: put $T_0':=0$ and let
\beaa
T_{k+1}'&:=& \min\{ n> T_k': S_{n-1} \in Conv(B), S_n \in B_+ \cup B_-\} \\
&& \wedge \min\{ n> T_k': S_{n-1} \in B_+, S_n \in B_+^c\} \, \wedge \, \min\{ n> T_k': S_{n-1} \in B_-, S_n \in B_-^c\}
\eeaa
for $k \ge 0$. Keep $H_k':=S_{T_k'}$ but put $\kappa':= \max\{k \ge 0: T_k' \le \tau_B \}$.
It is clear that if $B=(-d,d)$ and $|x| \ge d$, then $\kappa = \kappa'$ and $T_{k \wedge \kappa}=T_{k \wedge \kappa'}', H_{k \wedge \kappa}= H_{k \wedge \kappa'}'$ for any $k \ge 0$.

If $x \in B_+$, then $T_1'$ is the moment of the first entrance to the half-line $(-\infty, r)$ or $(-\infty, r]$ if $r \notin B$ or $r \in B$, respectively. The analogous is true for $x \in B_-$. Remarkably, \eqref{half-tail} is still valid with $T_1$ replaced by $T_1'$ regardless of whether the half-lines are open or closed; let us refer to this new equation as to (\ref{half-tail}$'$). Since this fact is not immediate, we discuss (\ref{half-tail}$'$) in Lemma~\ref{LEM UNIFORM'} in Sec.~\ref{SEC PROOFS MAIN} below.

Define
$$V_B(x):= \E_x \sum_{i =1}^\infty \bigl| S_{T_i'} - S_{T_{i-1}'} \bigr| \I_{\{T_i' \le \tau_B, \, S_{T_{i-1}'} \notin Conv(B) \}} = \E_x \sum_{i =1}^{\kappa'} \bigl| H_{i}' - H_{i-1}' \bigr| \I_{\{H_{i-1}' \notin Conv(B)\}},$$
where of course $V_{(-d,d)}(x) = V_d(x)$ and $V_B(x)\equiv 0$ on $B$.

\bt
\label{THM EQUIV GENERAL}
Let $S_n$ be a centred random walk with positive finite variance $\sigma^2$, and let $B \subset M$ a bounded Borel set such that $Int_M(B)$ is non-empty. Then for any $x \in M$,
\be \label{main B}
\P_x(\tau_B > n) \sim \frac{\sqrt{2} V_B(x)}{\sigma \sqrt{\pi n}}, \quad n \to \infty.
\ee
The function $V_B(x)$ is harmonic for the walk $S_n$ killed as it enters $B$, that is $V_B(x) = \E_x V_B(S_1)$ for $x \notin B$.
It holds that $0 \le V_B(x) < \infty$ for $x \notin B$ and $0< V_B(x)$ for $|x| \notin Conv(B)$, and $V_B(x) \sim |x|$ as $|x| \to \infty$. Moreover, \eqref{main B} is uniform in $x =o(\sqrt{n})$.
\et

It is unclear from the definition of $V_B(x)$ and \eqref{KS g=} that $V_B(x) \equiv \sigma^2 a_B(x)$ for integer-valued walks.

The following result describes the behaviour of the walk conditioned to avoid $B$ up to time $n$. Essentially, here we rigorously state the heuristics mentioned to comment on Theorem~\ref{THM EQUIV}. Define $\hat{S}_n(t)$ as follows: for $t=k/n$ with a positive integer $k$ put $\hat{S}_n(k/n):= S_k/(\sigma \sqrt{n})$, and define the other values by linear interpolation. Recall that for any $x \in B_+$ that is either fixed or $x=x_n$ such that $x_n = o(\sqrt{n})$,
\be \label{meander}
Law_x(\hat{S}_n(\cdot) | T_1'>n)  \stackrel{\D}{\to} W_+(\cdot) \qquad \mbox{in } C[0,1],
\ee
where $W_+$ is a standard Brownian meander. Let $\nu_n:=\max \{k \ge 0: T_k' \le n\}$ be the number of jumps over the edges of $B$ by the time $n$.

\bt
\label{THM MEANDER}
Under assumptions of Theorem~\ref{THM EQUIV GENERAL}, for any $x \in M$ such that $V_B(x) > 0$,
$$Law_x(\hat{S}_n(\cdot)| \tau_B>n) \stackrel{\D}{\to} \rho W_+(\cdot) \qquad \mbox{in } C[0,1],$$
where the sign $\rho$ of the meander is a random variable that is independent of $W_+$, with the distribution given by $\P(\rho = \pm 1) = \frac12 \pm \frac{x - \E_x S_{\tau_B}}{2 V_B(x)}$. Moreover, for any $b_n \to \infty$, we have that
\be \label{last jump 1}
\lim_{n \to \infty} \P_x(T'_{\nu_n} \le b_n | \tau_B>n ) = 1,
\ee
and for any $i \ge 1$,
\be \label{last jump 2}
\lim_{n \to \infty} \P_x({\nu_n} =i-1 | \tau_B>n ) = V_B^{-1}(x)\E_x \bigl| H_{i}' - H_{i-1}' \bigr| \I_{\{\kappa' \ge i, \, H_{i-1}' \notin Conv(B)\}}.
\ee

If $x=x_n$ and either $x_n \to \infty$ or  $x_n \to -\infty$ such that $x_n = o(\sqrt{n})$, then $\rho = 1$ or $\rho = -1$ a.s., respectively.
\et

For arithmetic random walks, this result is covered by Belkin~\cite{Belkin}, to be precise, by his Theorem 3.1, where the distribution of the limit process is stated incorrectly although the correct limit for one-dimensional distributions is given in Theorem 1.1. Again, his method can not be generalized for general walks. In the basic case that $B$ is an interval, the distribution of $\rho$ can be easily expressed in terms of $\E_x S_{\tau_B}$. This is very different in the general case, and we would not have provided a meaningful formula for $\P(\rho = 1)$ if we did not know Belkin's result.

\subsection{Applications to the largest gap problem}

Our original interest to the questions presented above was motivated by the problem of estimating size of the largest gap within the range of a random walk. Precisely, the largest gap $G_n$ by the time $n$ equals the maximal spacing, i.e.
$$G_n:= \max_{1 \le k \le n-1} S_{(k+1, n)} - S_{(k,n)},$$ where $m_n:=S_{(1, n)} \le S_{(2, n)} \le \dots \le S_{(n, n)}=:M_n$ denote the elements of $S_1, \dots, S_n$ arranged in the weakly ascending order. In particular, $G_n$ appears in the study of probability that the so-called iterated random walk $X(|S_k|)$ stays positive by the time $n$, where $X(t)$ is a L\'evy processes independent with the walk $S_k$. For this problem it is crucial to understand if the range of $S_n$ is dense enough to have it replaced by the continuous interval $[m_n, M_n]$. If so, the problem reduces to a much more feasible consideration of the probability that $X(t)$ stays above a negative level by the time $\max(|m_n|, M_n)$. We refer to Baumgarten~\cite{Baumgarten} and Vysotsky~\cite{MeIterated} for details.

Surprisingly, the author has found no references to the problem of estimating $G_n$. The closest result on sparseness of the range is of a rather different type and concerns the number of non-visited sites. It easily follows from Theorem 1.1 of Borodin~\cite{Borodin} (see his (2.1)) that for any $\lambda$-arithmetic centred random walk with finite variance, the quantity $$E_n:= \lambda^{-1}(M_n-m_n) - Card ( \{S_k \}_{k=1}^n) + 1,$$ which is the number of empty sites within the range, satisfies $$\frac{E_n}{\sqrt{n} } \stackrel{\P}{\longrightarrow} 0.$$

Clearly, the largest gap $G_n$ exceeds some value $h$ if and only if  by the time $n$ the walk has avoided a randomly located interval of length $h$. This idea and the $O(n^{-1/2})$ asymptotic from Theorem~\ref{THM EQUIV} are used in the argument proposed by Jian Ding and Yuval Peres to prove the following.

\begin{prop} \label{THM GAP}
If the random walk $S_n$ is centred and has a finite variance, then the family $Law(G_n)_{n \ge 1}$ is tight.
\end{prop}

Thus the largest gap is of the constant order. Actually, we can prove  more than this. Let $S_n^{\geqslant}$ and $S_n^<$ be {\it independent} Markov chains on $[0, \infty)$ and $(-\infty,0)$, respectively, defined by the transition probabilities $\P_x(S_1^{\geqslant} \in dy) = \frac{U_+(y)}{U_+(x)} \P_x(S_1 \in dy)$ for $x, y \ge 0$ and $\P_x(S_1^< \in dy) = \frac{U'_-(y)}{U'_-(x)} \P_x(S_1 \in dy)$ for $x \ge 0, y < 0$, where $U'_-(x) := \E_x S_{\tau_{[0, \infty)}} - x$. The Markov chain $S_n^{\geqslant}$ is a Doob $h$-transform of the random walk $S_n$ and equals the weak limit of the walk conditioned to stay positive: by Bertoin and Doney~\cite{BertoinDoney}, for any $k \ge 1$ and any bounded measurable $f:\R^k \to \R$,
\be \label{Doob h-limit}
\lim_{n \to \infty} \E_0 \bigl(f(S_1, \dots, S_k) \bigl| \bigr. \tau_{(-\infty, 0)}>n \bigr) = \E_0 f(S_1^{\geqslant}, \dots, S_k^{\geqslant}).
\ee
Analogously, the walk conditioned on $\tau'_{[0, \infty)}>n$ converges to $S_n^<$.

Let $G$ be a random variable distributed as the largest gap within the set $\{S_n^{\geqslant}, -S_n^<\}_{n \ge 0}$ under $S_0^{\geqslant} = S_0^<= 0$. Formally, this is the supremum of the spacings, which are well-defined since the chains $S_n^{\geqslant}$ and $-S_n^<$ are known to be transient, i.e. tend to infinity almost surely. If $S_n$ is $\lambda$-arithmetic with some $\lambda > 0$, let $E$ be a random variable distributed as the number of elements in $\lambda \N \setminus \{S_n^{\geqslant}, -S_n^<\}_{n \ge 0}$ under $S_0^{\geqslant} = S_0^<= 0$. So far it is not clear if $G$ and $E$ are proper random variables or they can be equal to infinity.

\bt
\label{THM GAP LIMIT}
Let $S_n$ be a centred random walk with a positive finite variance. For any sequence $b_n \to \infty$ such that $b_n=o(n)$ it holds that
\be \label{Int gap}
G_n^{Int}:= \max_{b_n \le k \le n-b_n} S_{(k+1, n)} - S_{(k,n)} \stackrel{\P}{\longrightarrow} u,
\ee
where $u = \lambda$ if the walk is $\lambda$-arithmetic and $u = 0 $ if the walk is non-arithmetic, and
\be \label{Ext gap}
G_n^{Ext}:= \max_{k \in [1, b_n] \cup [n-b_n, n-1] } S_{(k+1, n)} - S_{(k,n)} \stackrel{\D}{\longrightarrow} \max(G^-, G^+),
\ee
where $G^-$ and $G^+$ are i.i.d. positive proper random variables distributed as $G$. Moreover, if the walk is $\lambda$-arithmetic, then
\be \label{empty}
E_n \stackrel{\D}{\longrightarrow} E^- + E^+,
\ee
where $E^-$ and $E^+$ are i.i.d. proper random variables distributed as $E$.

\et
\begin{cor*}
$G_n \stackrel{\D}{\longrightarrow} \max(G^-, G^+)$.
\end{cor*}

Thus that the largest gap is attained near the extrema $m_n$ and $M_n$, and there are no gaps in the bulk. The later is related to the fact that asymptotically, a centred random walk has no points of increase, which itself is a discrete version a theorem by Dvoretzky, Erd{\"o}s and Kakutani that a Brownian motion never increases. A heuristic explanation of \eqref{Int gap} is given in Sec.~\ref{SEC PROOFS GAP} just before the proof of Theorem~\ref{THM GAP LIMIT}. Since an arithmetic walk visits all the points in the bulk, the non-visited sites are concentrated near the extrema, and the goal is to show that their total number is finite.

\medskip

The rest of the papers is organized as follows. The most important ideas and heuristics of our proofs of limit theorems on random walks that avoid a bounded set are given in Sec.~\ref{SEC IDEAS}. These theorems are proved in Sec.~\ref{SEC PROOFS MAIN}. The results on the largest gap are proved in Sec.~\ref{SEC PROOFS GAP}.

\section{Random walks that avoid a set: Ideas of proofs } \label{SEC IDEAS}

In this section we give an informal explanation of our proofs of Theorems~\ref{THM EQUIV}-\ref{THM MEANDER}. Let us start with Theorems~\ref{THM EQUIV} and \ref{THM MEANDER}. Consider the basic case. We want to argue that a typical strategy for the walk to avoid the interval $B=(-d,d)$ is to make few (possibly zero) jumps over $B$ within a short time and then stay on one side of $B$ until the end.

First, we claim that there exists a constant $\gamma (d) = \gamma \in (0,1)$ such that for every $|x| \ge d$,
\be \label{Jump over}
\P_x(|H_1| \ge d) \le \gamma.
\ee
This shows that each jump over $B$ contributes to probability of the trajectory a multiplicative factor that is bounded away from $1$.

In fact, $\E S_1^2 < \infty$ implies that $\E H_\pm < \infty$. Assume first that the walk in non-arithmetic. By the well-known result of renewal theory (Section XI.4 in Feller~\cite{Feller}), we get
\be \label{inf level}
\lim_{x \to -\infty} \P_x(|H_1| < d) = \frac{1}{\E H_+ } \int_0^{2d} \P( H_+  \ge t) dt >0,
\ee
where the r.h.s. corresponds to the distribution of the overshoot over an ``infinitely remote'' level. Thus \eqref{Jump over} holds true for all $x$ negatively large enough while for the ``small'' $x \le -d$ we simply force the walk to start with a sufficient (but uniformly bounded) number of negative jumps and thus reach a level that is remote enough. The case that $x \ge d$ should be considered in the same way. The $\lambda$-arithmetic case is completely analogous, with the only difference that
$$\lim_{x \to -\infty} \P_x(|H_1| < d) = \frac{1}{\E H_+ } \sum_{k=1}^{2 \lambda \lceil d/\lambda \rceil -1} k \P(H_+=k).$$
Hence \eqref{Jump over} is proved.

Since after each jump over the interval  the walk starts afresh, we need to control the size of the overshoot as seen from $p_n(x) \ge \P_x(T_1 >n) \ge c |x| n^{-1/2}$ (which follows by \eqref{half-tail}). For this purpose we will use the following estimate. For any $\alpha \in (0,1) $ there exists a constant $K(\alpha)$ such that for every $|x| \ge d$,
\be \label{Lyapunov func}
\E_x |H_1| \le \alpha |x| + K(\alpha).
\ee
For large $x$, this is by Theorem 3.10.2 of Gut~\cite{Gut}, which states that $\E H_1(x) = o(|x|)$ at infinity (and moreover, the family $\{H_1(x)/x \}_{|x| \ge d}$ is uniformly integrable, the fact to be used later). In order to use this reference we consider the random walk formed by the ladder heights of the walk $S_n$, with the integrable increments of the new walk distributed as $H_{sgn(x)}$.  For small $x$, \eqref{Lyapunov func} is immediate by $\E_x |H_1| \le U_{sgn(x)}(x - sgn(x)d)$ and the monotonicity of $U_\pm$.

Note that \eqref{Lyapunov func} actually gives a Lyapunov function for the Markov chain formed by the overshoots over the level $d$ from above, and the same is true for the chain over the overshoots of $-d$ from below.

Equation \eqref{Lyapunov func}, combined with \eqref{half-tail}, ensures that for the random walk that avoids the interval for a long time and at present is at a high level $x$, it is more efficient to stay at one side of the interval rather than jumping over and starting from the new level that on average is $o(x)$. In order to use this intuition  rigorously, we need to control $p_n(x)$ in terms of $x$. This is done in Lemma~\ref{LEM BOUND}, which states that $p_n(x) \le c |x| n^{-1/2}$.

Thus \eqref{Jump over} and \eqref{Lyapunov func} describe the mechanisms of contraction that lead to a typical absence of jumps over $B$ for large values of time on the event $\{\tau_B >n\}$. Having this idea understood, it is not hard to generalize Theorem~\ref{THM EQUIV} and prove Theorem~\ref{THM EQUIV GENERAL}. Essentially, the difference between $\tau_B$ and $\tau_{Conv(B)}$ stems from the time spent by the walk  by time $\tau_B$ in the ``holes'' inside $B$, i.e. $Conv(B) \setminus B$. Importantly, the exit time from any bounded set has an exponential tail, i.e. there exists constants $\beta, c>0$ such that for any random walk, for any $d >0$ and $|x| \le d$,
\be \label{two sided}
\P_x(\tau_{[-d,d]^c} > n) \le c e^{-\beta n};
\ee
so the time spend in the holes typically is negligible.

Equation \eqref{Lyapunov func} extends without any difficulties to its complete analogue (to be denoted as $(\ref{Lyapunov func}')$) that holds for any $x$. However, \eqref{Jump over} can not be generalized so easily as it is not sufficient to consider only one jump over the edges of $B$. This is readily seen from the following example: it holds that  $\P_0(\tau_{\{-1,1\}} > T_1') = 1$ if $X_1 \in \{-3,2\}$ a.s. We claim that for any bounded Borel set $B$ such that $Int_M(B) \neq \varnothing$ there exists a $\gamma' = \gamma'(B) \in (0,1)$ such that for any $x$,
\be \label{Jump over'}
\P_x(\tau_B > T_3') \le \gamma'.
\ee

It suffices to show existence of $\gamma'$ such that $\P_x(\tau_B > T_2') \le \gamma'$ for any $x \notin Conv(B)$ (because if $x \in Conv(B)$, simply exit from $Conv(B)$). In the arithmetic case this follows exactly as in the proof of \eqref{Jump over} because at each jump the walk must avoid $\{ l, r\} \subset B$. In the non-arithmetic case
a minor difficulty is that it may happen that the walk does not hit $B$ as it enters $Conv(B)$, e.g. if $B=(-1,1) \cup \{-2, 2\}$, $|X_1| \le 1$ a.s., and the distribution of $X_1$ is continuous. We resolve as follows: assume  w.l.o.g. that $x \in B_-$ and denote $q:= dist(Int(B), \{l\})$ and $h:= \sup \{s \in supp(X_1) \}$. Now choose an $\varepsilon \in (0, \lceil q/h \rceil h -q) $ such that $(l+q, l+q+\varepsilon) \in Int(B)$. Then for $k:= \lfloor q/h \rfloor $ it holds that $$\liminf_{x \to -\infty} \P_x(S_{T_1' + k} \in Int(B)) \ge \lim_{x \to -\infty} \P_x(S_{T_1'} \in (q - k h + \varepsilon, q - k h + 2\varepsilon)) \cdot  \P_0(S_k \in [ kh -\varepsilon, k h ]),$$ where the r.h.s. is strictly positive by \eqref{inf level} and $\sup \{s \in supp(H_+) \}=h$. Thus if starting from a distant $x$, the random walk hits $B$ after hitting $Conv(B)$ before exiting $Conv(B)$. For small $x$, argue as above for \eqref{Jump over}.

\section{Random walks that avoid a set: proofs of Theorems~\ref{THM EQUIV}-\ref{THM MEANDER}} \label{SEC PROOFS MAIN}

First note that  there exists a constant $c>0$ such that for any $n \ge $1 and $|x| \ge d$,
\be \label{Eq Eppel}
\P_x(T_1 > n) \le c |x| n^{-1/2}
\ee
(confer with \eqref{half-tail}). This follows, for example from Lemma 5 by Eppel~\cite{Eppel}, which gives a much stronger local version
\be \label{Eppel}
\P_x(T_1 = n) \le c \min(|x| n^{-3/2}, |x|^{-2}).
\ee

Now we prove the key estimate for $p_n(x)$.

\bl \label{LEM BOUND}
There exists a constant $c \ge 1$ such that $p_n(x) \le c |x| n^{-1/2}$ for any $x$ and $n$.
\el

\bp
The case that $|x| < d$ is trivial since $p_n(x)=0$ so assume that $|x| \ge d$. Fix an $\varepsilon \in (0,1)$ and let $n$ be large enough. Since $p_n(x)$ is monotone in $n$, write
$$p_n(x) \le \P_x(T_1 > \varepsilon n) + \E_x \, p_{(1- \varepsilon) n}(H_1) \I_{\{|H_1| \ge d\}}. $$ With the notation $I_x:= -sgn(x) [d, \infty)$, we have
\be \label{main bound}
p_n(x) \le \frac{C |x|}{\sqrt{\varepsilon n}} + \int_{I_x} p_{(1- \varepsilon) n}(y) \P_x (H_1 \in dy).
\ee
Now we estimate the integrand with \eqref{main bound} itself and use \eqref{Lyapunov func} to get that
\beaa
p_n(x) &\le& \frac{C |x|}{\sqrt{\varepsilon n}} + \int_{I_x} \left( \frac{C |y|}{\sqrt{\varepsilon (1- \varepsilon) n}} + \int_{I_y}   p_{(1- \varepsilon)^2 n}(z) \P_y (H_1 \in dz) \right) \P_x (H_1 \in dy)\\
&\le& \frac{C |x|}{\sqrt{\varepsilon n}} \Bigl( 1 + \frac{\alpha + K(\alpha)}{\sqrt{(1- \varepsilon)}}  \Bigr) + \int_{I_x} \int_{I_y} p_{(1- \varepsilon)^2 n}(z) \P_y (H_1 \in dz) \P_x (H_1 \in dy).
\eeaa
For the second iteration, use both  \eqref{Jump over} and \eqref{Lyapunov func} to obtain
$$p_n(x) \le \frac{C |x|}{\sqrt{\varepsilon n}} \Bigl( 1 + \frac{\alpha + K(\alpha)}{\sqrt{(1- \varepsilon)}} + \frac{\alpha^2 + \alpha K(\alpha) + \gamma K(\alpha)}{\sqrt{(1- \varepsilon)^2}} \Bigr) + \E_x p_{(1-\varepsilon)^3 n} (H_3) \I_{\{|H_i| \ge d, \, \, i=1,2,3 \}}$$

Letting $\varepsilon < 1 - \max(\alpha^2, \gamma^2)$, we repeat the estimation recursively. We can make at least $k$ recursions if $k$ satisfies $(1-\varepsilon)^k \varepsilon n - k \ge 2$, where the subtraction appears due to taking integer parts. Then $$p_n(x) \le \frac{C |x|}{\sqrt{\varepsilon n}} \Bigl( 1 - \frac{\alpha}{\sqrt{1-\varepsilon}} \Bigr)^{-1} +  \frac{CK(\alpha)}{\sqrt{(1-\varepsilon) n }} \Bigl( 1 - \frac{\gamma}{\sqrt{1-\varepsilon}} \Bigr)^{-1} \Bigl( 1 - \frac{\alpha}{\sqrt{1-\varepsilon}} \Bigr)^{-1} + \gamma^{k+1},$$
and by taking $\varepsilon $ to be sufficiently small and such that $\varepsilon < -2 \log \gamma$, we can choose $k$ to be large enough making the last term equal to $o(n^{-1/2})$.
\ep

We now give a strengthening of \eqref{half-tail}, which is needed to prove the uniform version of \eqref{main} as stated in Theorem~\ref{THM EQUIV}. Note that Lemma~\ref{LEM UNIFORM} is not required to prove \eqref{main} itself. The author thanks Ron Doney, who gave another proof of the missing part of the lemma, for the profound discussion on this problem.

\begin{lem} \label{LEM UNIFORM}
For any centred random walk with a positive finite variance $\sigma^2$, for any $d >0$ \eqref{half-tail} holds uniformly in $x=o(\sqrt{n})$.
\end{lem}

\bp
For the walks that are not non-centred lattice, this follows immediately from Theorem 2 in Doney~\cite{Doney}, which gives the uniform asymptotic for the local probabilities $\P_x( T_1 = n)$. Recall that the distribution of $X_1$ is $(h, a)$-lattice with  the span $h>0$ if $h$ is the maximal number such that $\P(X_1 \in b + h \Z) = 1$ for some $b \in \R$ and the shift $a \in [0, h)$ is such that $\P(X_1 \in a + h \Z) = 1$; the lattice is non-centred if $a \neq 0$.

The non-centred lattice case is very different as the local asymptotic does not match the tail asymptotic: Theorem 20 of~\cite{Doney} states that there exists a non-negative function $D(x)$ on $[0, h)$ such that
\be \label{local Doney}
\P_{-x}(\tau_+ =  n+1 ) \sim \frac{\sqrt{2} U_-(x) D( (n (h-a) +x)^*)}{2\sigma \sqrt{\pi} n^{3/2}}
\ee
uniformly in $x = o(\sqrt{n})$, where $y^*:=y - h \lfloor y/h \rfloor$ denotes the $h$-fractional part of $y \ge 0$. Then the required statement easily follows by standard summation once we check that
\be \label{averaging}
\lim_{n \to \infty} \sup_{0 \le y \le h } \Bigl | \frac{1}{n}  \sum_{k=1}^n D( (k (h-a) +y)^*) - \bar{D}(y) \Bigr | = 0
\ee
for some function $\bar{D}(y)$; then necessarily $\bar{D}(y) \equiv 1$ since we already know the coefficient in \eqref{half-tail}. Under the supremum we made the change $y= (l(h-a)+x)^*$, where $0 \le x \le h, l \ge n$.

By (74) of~\cite{Doney} we have that
\be \label{D}
D(x)= \sum_{m=0}^\infty V(x+hm) \bar{F}(x+hm), \quad x \in [0, h],
\ee
where $\bar{F}(y):= \P(X_1 >y)$ and $V(y)$ is the renewal function for weak descending ladder height process. Since both functions are monotone, (73) of~\cite{Doney} implies that first, $D(x)$ is bounded away from infinity, and second, by the Weierstrass' criterion, the series in \eqref{D} converges uniformly. Importantly, $V(x)$ and $\bar{F}(x)$ are monotone and right-continuous. Then $D(x)$ is a cadlag function as a uniform limit of cadlag functions by the completeness of the Skorokhod space under the metric $d_0$, see Sec. 12 of Billingsley~\cite{Billingsley}.

If $a/h$ is rational and equals $p/q$ for some positive integers $p, q$ (so the distribution is $h/q$-arithmetic), then \eqref{averaging} clearly holds true with $\bar{D}(x) = q^{-1}\sum_{k=0}^{q-1} D((k/q + x)^*)$ since $D(x)$ is bounded.

For an irrational $a/h$, a classical result of H. Weyl states that the sequence $(n(1-a)+y)^*$ is uniformly distributed on $[0, h]$ and that for any Riemann integrable $D(y)$, the averages in \eqref{averaging} converge to $\bar{D} (y) \equiv h^{-1} \int_0^h D(x) dx$ for any fixed $y$; see, e.g. Problems 162-166 in Part II, Sec. 4 of Polya and Szego~\cite{Problems}. We only need to check the uniformity of this convergence. By Lemma 1 in Sec. 12 of Billingsley~\cite{Billingsley}, for any cadlag function $D(x)$ there exists two step functions (finite linear combinations of indicator functions of intervals) $D_\pm^{(\varepsilon)}(x)$ such that $D_-^{(\varepsilon)} (x) \le D(x) \le D_+^{(\varepsilon)} (x)$ on $[0,h]$ and $D_+^{(\varepsilon)} (x) - D_-^{(\varepsilon)} (x) \le \varepsilon$. Hence it suffices to prove \eqref{averaging} for any step function. This easily follows from the fact that the empirical distribution functions $\# \{1 \le k \le n: (k(1-a))^* \le x \}/n$ converge to $x/h$ not only point-wisely but uniformly in $x \in [0,h]$.

\ep

\bp[{\bf Proof of Theorem~\ref{THM EQUIV}}]

Let us first show that
\be \label{1 immediate}
p_n(x) \sim \frac{\sqrt{2} \E_x |x - H_1|}{\sigma \sqrt{\pi n} } + \E_x p_n(H_1) \I_{\{ |H_1| \ge d\}}
\ee
uniformly in $d \le |x| \le a_n \sqrt{n}$.

We claim that for any sequence $b_n$ such that $b_n=o(n)$ and $a_n^2 n = o(b_n)$ (note that $b_n \to \infty$),
\be \label{total pb}
p_n(x) \sim \P_x(T_1 > n) + \P_x( T_1 \le b_n, \tau > n).
\ee
holds uniformly in $d \le |x| \le a_n \sqrt{n}$. Clearly, it suffices to check that
\be \label{not immediate}
\P_x( b_n < T_1 \le n, \tau > n) = o \Bigl(\frac{|x|}{\sqrt{n}} \Bigr)
\ee
uniformly. This means that conditioned on $\{\tau > n \}$, the walk either jumps over the interval ``immediately'' or stays to the same side of $(-d,d)$ for the whole time.

For any $\delta \in (0,1/2)$,
\beaa
\P_x( b_n < T_1 \le n, \tau > n)  &\le& \P_x( (1-\delta)  n < T_1 \le n) + \P_x(   b_n < T_1 \le (1- \delta) n, \tau > n) \\
&\le& 2C \Bigl (\frac{1}{\sqrt{1-\delta}} -1 \Bigr) \frac{|x|}{\sqrt{n}} + \P_x(   b_n < T_1 \le (1- \delta) n, \tau > n),
\eeaa
where we used \eqref{Eppel}. The contribution of the first term vanishes as  $\delta \to 0+$. For the second, by Lemma~\ref{LEM BOUND}, it holds that
\bea
\P_x(   b_n < T_1 \le (1- \delta) n, \tau > n) &\le& \int_{I_x} p_{\delta n} (y) \P_x(H_1 \in dy, T_1 >  b_n ) \label{est 1}\\
&\le& \frac{c \E_x |H_1| \I_{\{ T_1 >  b_n \}}}{\sqrt{\delta n}} \notag \\
&\le& \frac{c|x|} {\sqrt{\delta n} }  \E \left[ \frac{|H_1(x)|}{|x|} \I_{ \{T_1 (sgn(x) a_n \sqrt{n}) >  b_n \} } \right], \notag
\eea
where the indicator tends to $0$ a.s. Then the expectation converges to zero as $n \to \infty$ uniformly in $x$ by the uniform integrability of $\{ H_1(x)/x \}_{|x| \ge d}$. Thus \eqref{not immediate} is true and \eqref{total pb} follows.

By Lemma~\ref{LEM UNIFORM}, the first term in the r.h.s. of \eqref{total pb} gives the corresponding one in \eqref{1 immediate}, and it remains to consider the second term. It is upper bounded by $\E_x p_{n-b_n} (H_1) \I_{\{|H_1| \ge d \}}$
while for a lower bound,
\beaa
\P_x( T_1 \le b_n, \tau > n)  &\ge& \E_x p_n (H_1) \I_{\{|H_1| \ge d, \, T_1 \le b_n \}} \\
&=& \E_x p_n (H_1) \I_{\{|H_1| \ge d \}} - \E_x p_n (H_1) \I_{\{|H_1| \ge d, \, T_1 > b_n \}},
\eeaa
where the subtracted term is dominated by the r.h.s. of \eqref{est 1} and thus its contribution vanishes as $n \to \infty$. We estimate the difference of the main terms using \eqref{Eppel} and get
\bea
\E_x \bigl( p_{ n - b_n} (H_1) -  p_{n} (H_1) \bigr) \I_{\{|H_1| \ge d \}} &\le& \int_{I_x} \P_y(n - b_n < T_1 < n) \P_x(H_1 \in dy) \label{difference}\\
&\le& \int_{I_x} \sum_{i=n-b_n +1}^n \frac{C |y|}{i^{-3/2}} \P_x(H_1 \in dy) \notag \\
&\le& \frac{3 C \E_x |H_1|}{|x|} \cdot \frac{b_n}{n} \cdot \frac{|x|}{\sqrt{n}} \notag
\eea
for $n$ large enough, hence \eqref{1 immediate} follows from \eqref{total pb}.

We now claim that
\be \label{2 immediate}
\E_x p_n(H_1) \I_{\{ |H_1| \ge d\}} \sim \frac{\sqrt{2} \E_x |H_2 - H_1| \I_{\{ |H_1| \ge d\}}  }{\sigma \sqrt{\pi n} } + \E_x p_n(H_2) \I_{\{ |H_i| \ge d, \,\, i=1,2 \}}
\ee
uniformly over $d \le |x| \le a_n \sqrt{n}$. Indeed, \eqref{1 immediate} can be applied on the set $\{ d \le |H_1| \le a_n \sqrt{n} \}$ and \eqref{2 immediate} follows if we show that the contribution of the expectations over
the complement set vanishes as $n \to \infty$. We have that
$$\E_x p_n(H_1) \I_{\{ |H_1| > a_n \sqrt{n} \}} \le \frac{c \E_x |H_1| \I_{\{ |H_1| > a_n \sqrt{n} \}} }{\sqrt{n}} \le \frac{c |x|}{\sqrt{n}} \E \left [\frac{|H_1(x)|}{|x|} \I_{ \bigl \{ \frac{|H_1(x)|}{|x|} > \frac{a_n \sqrt{n}}{|x|} \bigr \}} \right],$$ and the expectation converges to $0$ if $x=o(a_n \sqrt{n})$ by the uniform integrability while in the zone $|x| \to \infty$ we simply use that $\E_x|H_1| =o(|x|)$. Since
$$\E_x |H_2 | \I_{\{ |H_1| \ge a_n \sqrt{n}\}} \le \E_x c_n |H_1| \I_{\{ |H_1| \ge a_n \sqrt{n}\}}$$ for some $c_n \to 0$ as $n \to \infty$, the rest follows. Note that it suffices to use \eqref{half-tail} rather than Lemma~\ref{LEM UNIFORM} to prove \eqref{2 immediate} for a fixed $x$ as follows from Lemma~\ref{LEM BOUND} and the dominated convergence theorem.

We apply \eqref{1 immediate} to \eqref{2 immediate}, and do so recursively to obtain that
\be \label{k immediate}
p_n(x) \sim \frac{\sqrt{2} \E_x \sum_{i=1}^{\min(k,  \kappa)} \bigl| H_{i} - H_{i-1} \bigr|}{\sigma \sqrt{\pi n}} + \E_x p_n(H_k)   \I_{\{\kappa \ge k +1 \}}
\ee
holds uniformly for $x = o(\sqrt{n})$ for every fixed $k \ge 1$. It remains to note that
\be \label{tail 1}
\E_x |H_i| \I_{\{ \kappa \ge i \}} \le \alpha^i |x| + K(\alpha) (\alpha^{i-1} + \alpha^{i-2} \gamma + \dots + \gamma^{i-1}),
\ee
which is obtained by a recursive application of \eqref{Jump over} and \eqref{Lyapunov func}. Similarly, by Lemma~\ref{LEM BOUND}
\be \label{tail 2}
\E_x p_n(H_k) \I_{\{\kappa \ge k +1 \}} \le \frac{C}{\sqrt{n}} \bigl(\alpha^k |x| + K(\alpha) (\alpha^{k-1} + \alpha^{k-2} \gamma + \dots + \gamma^{k-1}) \bigr),
\ee
thus this term vanishes as $k \to \infty$. A simple computation then gives
$$\lim_{k \to \infty} \E_x \sum_{i=1}^{\min(k,  \kappa)} \bigl| H_{i} - H_{i-1} \bigr| = V(x) \le \frac{1+ \alpha}{1-\alpha} |x| + \frac{2K(\alpha)}{({1-\alpha} )(1-\gamma)},$$ which is particular shows that $V(x)$ is finite and $V(x) \sim |x|$ at infinity. The estimates above and \eqref{k immediate} conclude the proof of the uniform version of \eqref{main}.

Further, since $V_d(\pm(d + y)) = V_{\pm (-2d, 0)}(y)$ for any $d >0$ and any fixed $y \ge 0$, we have that $$V_d(\pm(d + y)) - U_d(\pm (d +y)) = \E_{\pm y} \sum_{i=2}^{\kappa'} \bigl| H'_{i} - H'_{i-1} \bigr|,$$ where the random variables $\kappa'$ and $H'_i$ are considered for the set $B = \pm (-2d, 0)$. The integrand monotonously converges to $0$ a.s. as $d \to \infty$ because $\kappa' \to 1$ for almost every trajectory of the walk, and hence the expectation tends to zero.

We will prove the harmonicity of $V_{(-d,d)}(x)$ in Theorem~\ref{THM EQUIV GENERAL}.
\ep

Recall that Theorem~\ref{THM EQUIV GENERAL}, which generalizes Theorem~\ref{THM EQUIV}, is stated for a set $B$ regardless whether it includes or does not include its boundary points $l$ and $r$. As we mentioned in the introduction, \eqref{half-tail} is still valid with $T_1$ replaced by $T_1'$. Since this is not straightforward and no reference is available, we state this relation separately:

\begin{lem} \label{LEM UNIFORM'}
Let $S_n$ be a centred random walk with a positive finite variance $\sigma^2$, and let $B$ be a non-empty bounded Borel set. For any fixed $x \in M \setminus Conv(B)$ we have that
$$\P_x( T_1' >n) \sim \frac{\sqrt{2} |x - \E_x S_{T_1'}|}{\sigma \sqrt{\pi n} }, \quad n \to \infty. \eqno{(\ref{half-tail}')}$$
Moreover, this holds uniformly in $x=o(\sqrt{n})$.
\end{lem}
\bp
It suffices to consider only $B=\{0\}$ since the case that the boundary points are not included in $B$ is already covered by \eqref{half-tail} and Lemma~\ref{LEM UNIFORM}. Denote $H_\pm':= S_{\tau_{\pm [0, \infty)}}$, $U'(x):=|x - \E_x S_{T_1'}|$ for $x \neq 0$, $U_+'(x):= x - \E_x H_-'$ for $x \ge 0$, $U_-'(x):= \E_x H_+' - x$ for $x \le 0$. W.l.o.g. assume that the walk is non-arithmetic (otherwise simply use that $\tau_{\{0\}} = \tau_{(-\lambda, \lambda)}$ for $x \in M$).

Again, by Feller~\cite{Feller}, it holds that $\P_0(\tilde{\tau}_{\pm [0, \infty)} > n) \sim \frac{\sqrt{2}}{\sigma \sqrt{\pi n} } U_\pm'(0)$ so the coefficients for positive and negative $x$ match, and we can combine the tail asymptotics for the lengths of negative and positive excursions into a single formula; here we used $\tilde{\tau}_{\pm [0, \infty)} $ rather than $\tau_{\pm [0, \infty)}$ since the latter equals $0$ under $\P_0$. Now w.l.o.g., consider $x <0$. We claim that $U_-'(x) = U_-(x+)$ for $x<0$. We already know that \eqref{half-tail} is true if $x$ is a point of continuity of $U_-(x)$, so assume that there is a discontinuity at $x$.

Clearly, it suffices to prove that
\be \label{diff sim 1}
\P_x(\tau_{(0, \infty)} > n) - \P_x(\tau_{[0, \infty)} > n) \sim \frac{\sqrt{2} (U_-(x) - U_-(x+))}{\sigma \sqrt{\pi n} }.
\ee
Let $H_k$ be the random walk of the strong ascending ladder heights, i.e. a walk with the increments distributed as $H_+$, and put $H_0:=0$. By the Wald identity, we have that $U_-(x) = \E H_+ \sum_{k=0}^\infty \P(H_k \le -x)$ for $x \le 0$, where the sum equals the renewal function. Then  \eqref{diff sim 1} is equivalent to
\be \label{diff sim 2}
\lim_{n \to \infty} \sum_{i=1}^n \P_x(T_1' = i, S_{T_1'} = 0) \Bigl ( \sqrt{n} \P_0(\tau_{(0,\infty)} > n-i) \Bigr) = \frac{\sqrt{2} \E H_+}{\sigma \sqrt{\pi}} \sum_{k=1}^\infty \P(H_k = -x).
\ee

Let us split the l.h.s. in two parts $\Sigma_1$ and $\Sigma_2$ that correspond to $i \le n/2$ and $i > n/2$, respectively. The summands of $\Sigma_1$ are dominated term-wisely by the summable sequence $c \P_x(T_1' = i, S_{T_1'} = 0)$ for some $c>0$. Since $\sum_{i=1}^\infty \P_x(T_1' = i, S_{T_1'} = 0) = \sum_{k=1}^\infty \P(H_k = -x)$, the dominated convergence theorem and \eqref{half-tail} imply that $\Sigma_1$ converges to the r.h.s. of \eqref{diff sim 2}.

For the remaining term, for any $R<0$ we have that
\beaa
\Sigma_2 &\le& \sqrt{n} \P_x(S_{T_1'} = 0, T_1' > n/2) \\
&\le& \sqrt{n} \P_x(T_1' > n/2) \Bigl [ \int_{-\infty}^{R} \P_y(S_{T_1'} = 0) \P_x(S_{n/2} \in dy | T_1' > n/2) + \P_x(S_{n/2} >R | T_1' > n/2) \Bigr ].
\eeaa
For any fixed $x$, the factor in front of the brackets is bounded. The second term in the brackets tends to zero for any fixed $R<0$. The first term can be made as small as necessary by choosing $R$ to be negatively large enough. Indeed, $\lim_{y \to -\infty} \P_y(S_{T_1'} = 0) = 0$ as the overshoot of a non-arithmetic walk over an infinitely remote level has a continuous distribution, cf. \eqref{inf level}. Thus $(\ref{half-tail}')$ is proved.

Now check that the convergence is uniform. First note that Lemma~\ref{LEM UNIFORM}, \eqref{diff sim 1} and the monotonicity of $\P_x( T_1' >n)$ in $x$ imply that $(\ref{half-tail}')$  holds uniformly over any compact set. The rest follows by  Lemma~\ref{LEM UNIFORM} and the fact that for non-arithmetic walks, $U_-(x)$ asymptotically has no jumps, in the sense that $\lim \limits_{R \to \infty} \sup \limits_{x \le -R} (U_-(x)-U_-(x+)) = 0$. Again, this is true by the absence of atoms in the distribution of the overshoot over an infinitely remote level.
\ep

\bp[{\bf Proof of Theorem~\ref{THM EQUIV GENERAL}}]

The proof closely follows the one of Theorem~\ref{THM EQUIV}. W.l.o.g., assume that $0 \in Int(B)$. For $p_n'(x):= \P_x(\tau_B >n)$, we show that
\be \label{1 immediate'}
p_n'(x) \sim \frac{\sqrt{2} \E_x |x - H_1'| \I_{\{x \notin Conv(B) \}}}{\sigma \sqrt{\pi n} } + \E_x p_n'(H_1') \I_{\{ \tau_B > T_1'\}}
\ee
uniformly in $x \notin B$ and $|x| \le a_n \sqrt{n}$, which is analogous to \eqref{1 immediate}. Indeed, if $x \in Conv(B)$, then \eqref{two sided} easily implies that $p_n'(x) \sim \E_x p_n'(H_1') \I_{\{ \tau_B > T_1'\}}$, otherwise we proceed using Lemma~\ref{LEM UNIFORM'} exactly as in \eqref{1 immediate}. The difference is that we use $(\ref{Lyapunov func}')$ instead of \eqref{Lyapunov func} and get the estimate $p_n'(x) \le C |x|/\sqrt{n}$ by applying Lemma~\ref{LEM BOUND} to an interval that is contained in $Int(B)$.

Then we apply \eqref{1 immediate} recursively and see that
$$p_n'(x) \sim \frac{\sqrt{2} \E_x \sum_{i=1}^{\min(k,  \kappa')} \bigl| H_{i}' - H_{i-1}' \bigr| \I_{\{H_{i-1}' \notin Conv(B) \}}}{\sigma \sqrt{\pi n}} + \E_x p_n'(H_k') \I_{\{\tau_B > T_k' \}}$$
holds uniformly in $x = o(\sqrt{n})$ for any $k\ge 1$. By a recursive application of $(\ref{Lyapunov func} ')$ and \eqref{Jump over'}, we obtain
\be
\label{tail 1'}
\E_x |H_{3i}'| \I_{\{ \tau_B > T_{3i}' \}} \le \alpha^i |x| + K'(\alpha) (\alpha^{i-1} + \alpha^{i-2} \gamma' + \dots + (\gamma')^{i-1}),
\ee
which is analogous to \eqref{tail 1}, while
\be
\label{tail 2'}
\E_x |H_{3i+ j}'| \I_{\{ \tau_B > T_{3i+j}' \}} \le 2 K'(\alpha) + \alpha \E_x |H_{3i}'| \I_{\{ \tau_B > T_{3i}' \}}, \quad j=1,2.
\ee
Then we apply the estimate $p_n'(x) \le C |x|/\sqrt{n}$ to get an analogue of \eqref{tail 2}, and the rest follows.

All the stated properties of $V_B(x)$ except the harmonicity are straightforward. To prove the latter, we note that for any $x \notin B$ and $y \in \R$,
$$\P_x \bigl( S_{T_1' \wedge \tau_B}, S_{T_2' \wedge \tau_B}, \dots \bigl | \bigr. S_1 = y \bigr) = \P_y \bigl( S_{T_1' \wedge \tau_B}, S_{T_2' \wedge \tau_B}, \dots \bigr),$$ which simply means that any trajectory of the walk shifted one step forward enters the fixed sets $B_+$, $B_-$ and $Conv(B)$ at the same locations as the original trajectory.

The relation $V_B(x) =\E_x V_B(S_1)$ for $x \in B_+$ follows from
\beaa
&& V_B(x) \\
&=& \E_x \sum_{i=1}^\infty |S_{T_i' \wedge \tau_B} - S_{T_{i-1}' \wedge \tau_B}| \I_{\{ S_{T_{i-1}' \wedge \tau_B} \notin Conv(B) \}} \\
&=& \int_{-\infty}^\infty \E_x \Bigl [\sum_{i=2}^\infty |S_{T_i' \wedge \tau_B} - S_{T_{i-1}' \wedge \tau_B}| \I_{\{ S_{T_{i-1}' \wedge \tau_B} \notin Conv(B) \}} + x -S_{T_1'  \wedge \tau_B } \Bigl | \Bigr. S_1 = y \Bigr] \P_x(S_1 \in dy ) \\
&=& \int \limits_{-\infty}^\infty x- y + \E_x \Bigl [\sum_{i=2}^\infty |S_{T_i' \wedge \tau_B} - S_{T_{i-1}' \wedge \tau_B}| \I_{\{ S_{T_{i-1}' \wedge \tau_B} \notin Conv(B) \}} + y -S_{T_1'  \wedge \tau_B} \Bigl | \Bigr. S_1 = y \Bigr] \P_x(S_1 \in dy ) \\
&=& \int_{-\infty}^\infty (x- y + V_B(y)) \P_x(S_1 \in dy ),
\eeaa
where the last equation is valid because first, both the conditional expectation and $V_B(y)$ equal $0$ for $y \in B$ and second, $S_1 \ge S_{T_{i-1}' \wedge \tau_B}$ a.s. under $\P_x$ with $x \in B_+$ (since $S_{T_{i-1}' \wedge \tau_B} \le r$ suffices if $S_1 \in B_+$ and $S_1 = S_{T_{i-1}' \wedge \tau_B}$ if $S_1 \notin B_+$). The case that $x \in B_-$ is analogous. The remaining case that $x \in Conv(B) \setminus B$ is even simpler as the term in the sum for $i=1$ equals zero.

\ep

\bp[{\bf Proof of Theorem~\ref{THM MEANDER}}]
First consider the case of a fixed $x$. W.l.o.g., assume that $0 \in Int(B)$ and $b_n=o(n)$. As usual, we need to check convergence of finite-dimensional conditional distributions and tightness.

\underline{Convergence of finite-dimensional distributions.} The idea is straightforward so we check only for one-dimensional distributions to avoid bulky notation. We claim that for any $u \ge 0$ and $0 <t \le 1$,  $p_n''(x):= \P_x(0 \le \hat{S}_n(t) \le u, \tau_B > n)$ satisfies
\be \label{1D}
\lim_{n \to \infty} \sqrt{n} p_n''(x) = \frac{\sqrt{2} \P( W_+(t) \le u)}{\sigma \sqrt{\pi }} \E_x \sum_{i =1}^{\kappa'} \bigl| H_{i}' - H_{i-1}' \bigr| \I_{\{H_{i-1}' \in B_+\}}.
\ee

Fix a $k \ge 1$ and define $\ell_n:= \min\{j \ge 1: T_j' > j b_n/ (\nu_n \wedge k)\}$ if $\nu_n \ge 1$ and $\ell_n:=0$ if $\nu_n=0$. Clearly, $\ell_n \le \nu_n \wedge k$ on the event $\{T'_{\nu_n \wedge k} > b_n \}$, hence
\beaa
\Bigl \{ T'_{\nu_n \wedge k} > b_n, \tau_B > n \Bigr \} &=& \bigcup_{j=1}^k \Bigl \{  \ell_n = j, T'_{\nu_n \wedge k} > b_n , \tau_B > n \Bigr \} \\
&\subset& \bigcup_{j=1}^k \Bigl \{ T_j' - T_{j-1}' > \frac{b_n}{\nu_n \wedge k}, T_{j-1}' \le \frac{j-1}{\nu_n \wedge k} b_n, j = \ell_n \le \nu_n, \tau_B > n \Bigr \} \\
&\subset& \bigcup_{j=1}^k \Bigl \{ T_j' - T_{j-1}' > b_n /k , T_{j-1}' \le b_n, T_j' \le n < \tau_B \Bigr \}.
\eeaa

Now argue by analogy with the proof of \eqref{not immediate}: for any $\delta \in (0, 1/2)$,
\beaa
\Bigl \{ T'_{\nu_n \wedge k} > b_n, \tau_B > n \Bigr \} &\le& \bigcup_{j=1}^k \Bigl \{ (1-\delta) n \le T_j' - T_{j-1}' \le n, \tau_B > T_{j-1}' \Bigr \} \\
&& \cup  \, \Bigl \{ T_j' - T_{j-1}' > b_n/k, T_j' < (1-\delta) n + b_n, \tau_B > n \Bigr \}
\eeaa
and by conditioning on $H_{j-1}', T_{j-1}'$ for the events in the r.h.s. of the first line and on $H_{j-1}', T_{j-1}', T_j'$ for the events in the second line, we obtain
\beaa
\P_x \bigl ( T'_{\nu_n \wedge k} > b_n, \tau_B > n  \bigr ) &\le& \frac{2 C}{\sqrt{n}} \Bigl ( \frac{1}{\sqrt{1-\delta}} - 1\Bigr) \E_x \sum_{j=1}^\infty |H_{j-1}'| \I_{\{\tau_B > T_{j-1}'\}}  \\
&+& \frac{c}{\sqrt{\delta n}} \E_x \sum_{j=1}^\infty |H_{j-1}'| \I_{\{\tau_B > T_{j-1}'\}} \I_{\{ \max \limits_{1 \le i \le k} T_i' - T_{i-1}' > b_n/k\}}.
\eeaa

The expectation of the sum is finite by \eqref{tail 1'} and \eqref{tail 2'}, and the second indicator converges to $0$ a.s. for any fixed $k$, hence
$$\lim_{n \to \infty} \sqrt{n} \P_x \bigl ( T'_{\nu_n \wedge k} > b_n, \tau_B > n  \bigr ) = 0.$$
Further, $$\P_x \bigl ( T'_{\nu_n \wedge k} \le b_n, \nu_n \ge k+1, \tau_B > n  \bigr ) \le \P_x \bigl ( T'_k \le b_n, \tau_B > n  \bigr ) = \E_x p'_{n - b_n}(H_k') \I_{\{ \tau_B > T_k'\}},$$
Lemma~\ref{LEM BOUND} (recall that $B$ contains a centred interval) and $\lim \limits_{k \to \infty} \E_x |H_k'| \I_{\{ \tau_B > T_k'\}} = 0$, which follows by \eqref{tail 1'} and \eqref{tail 2'}, imply that for any events $A_n$ it holds that
\be \label{truncation}
\lim_{n \to \infty} \P_x \bigl ( A_n | \tau_B > n  \bigr ) = \lim_{k \to \infty} \lim_{n \to \infty} \P_x \bigl ( A_n, T'_{\nu_n} \le b_n, \nu_n \le k | \tau_B > n  \bigr )
\ee
given that the limits in $n$ exist. Now \eqref{last jump 1} follows immediately, and we also see that the family $\{ Law_x(\nu_n | \tau_B >n) \}_{n \ge 1}$ is tight.

In order to prove \eqref{1D}, we study the behaviour of $$\P_x(0 \le \hat{S}_n(t) \le u, T'_{\nu_n} \le b_n, \nu_n = i-1, \tau_B >n)$$ for $i \ge 1$. Let us
simplify this expression as we did in the proof of Theorem~\ref{THM EQUIV}. The main difference is that more effort is needed to get the analogue of \eqref{difference}. For any sequence $a_n$ such that $a_n  \to 0$ and $a_n \sqrt{n} \to \infty$ and any fixed $i \ge 1$, it holds that
$$\P_x(|H'_{i-1}| > a_n \sqrt{n}, T'_{i-1} \le b_n, \tau_B >n) \le \E_x p'_{n-b_n}(H'_{i-1}) \I_{\{|H'_{i-1}| > a_n \sqrt{n}\}}= o_i(n^{-1/2}).$$ By \eqref{two sided}, we have
\beaa
&& \P_x(0 \le \hat{S}_n(t) \le u, T'_{\nu_n} \le b_n, \nu_n = i-1, \tau_B >n)  \\
& \sim & \sum_{j=0}^{b_n} \int_r^{a_n \sqrt{n}}  \P_x(0 \le \hat{S}_n(t) \le u, T'_{i-1} = j, H'_{i-1} \in dy, T'_i >n) \\
& = & \sum_{j=0}^{b_n} \int_r^{a_n \sqrt{n}}  \P_y \Bigl(0 \le \frac{S_{nt - j}}{\sigma \sqrt{n}} \le u, T'_1 >n-j \Bigr)
\P_x(\tau_B > T'_{i-1} = j, H'_{i-1} \in dy,) \\
&\le& \sum_{j=0}^{b_n} \int_r^{a_n \sqrt{n}} \P_y \Bigl( \min_{nt - b_n \le \ell \le nt} \frac{S_\ell}{\sigma \sqrt{n}} \le u , T'_1 >n-b_n\Bigr) \P_x(\tau_B > T'_{i-1} = j, H'_{i-1} \in dy),
\eeaa
for an upper bound, and the analogous expression for a lower bound with the integrand replaced by $\P_y \bigl( \max \limits_{nt - b_n \le \ell \le nt} \frac{S_\ell}{\sigma \sqrt{n}} \le u, T'_1 > n \bigr) $.

Now by \eqref{meander}, which holds uniformly in the given range of $y$, and Lemma~\ref{LEM UNIFORM'} we get
\beaa
&& \P_x(0 \le \hat{S}_n(t) \le u, T'_{\nu_n} \le b_n, \nu_n = i-1, \tau_B >n)  \\
&\sim & \frac{\sqrt{2} \P ( W_+(t) \le u )}{\sigma\sqrt{\pi n}} \E_x (H_i' -H'_{i-1}) \I_{\{\kappa' \ge i, H'_{i-1} \in B_+ \}} \I_{\{H'_{i-1} \le a_n \sqrt{n}, T'_{i-1} \le b_n \}}á
\eeaa
where the second indicator tends to $1$ a.s. It remains to sum over $i$ from $1$ to $k$. Then \eqref{1D} follows by \eqref{truncation}, while \eqref{last jump 2} follows by adding up the remaining trajectories with negative $\hat{S}_n(t)$.

\underline{Tightness.} It suffices (see Billingsley~\cite[Sec. 7]{Billingsley}) to check that for any $\varepsilon >0$,
\be \label{tightness}
\lim_{\delta \to 0} \liminf_n \P_x( \omega_\delta(\hat{S}_n) \le \varepsilon| \tau_B >n) = 1,
\ee
where $\omega_\cdot(f)$ is the modulus of continuity of a function $f$ on $[0,1]$. Analogously to the proof of \eqref{1D}, for any $i \ge 1$ and $n$ large enough,
\beaa
&& \P_x( \omega_\delta(\hat{S}_n) \le \varepsilon, T'_{i-1} \le b_n, \tau_B \ge T'_i >n)  \\
& \ge & \sum_{j=0}^{b_n} \int_r^{a_n \sqrt{n}} \P_x( \omega_\delta(\hat{S}_n) \le \varepsilon, T'_{i-1} = j, H'_{i-1} \in dy, T'_i >n)\\
& \ge & \sum_{j=0}^{b_n} \int \limits_r^{a_n \sqrt{n}}  \P_y \Bigl( \max_{0 \le \ell < m \le n-j \atop m-\ell \le \delta n} \frac{|S_m - S_\ell |}{\sigma \sqrt{n}} \le \frac{\varepsilon}{2}, T'_1 >n-j \Bigr) \P_x \Bigl(\max_{0 \le \ell \le j} \frac{|S_\ell|}{\sigma \sqrt{n}} \le \frac{\varepsilon}{4}, \tau_B > T'_{i-1} = j, H'_{i-1} \in dy \Bigr) \\
& \ge & \int_r^{a_n \sqrt{n}}  \P_y ( \omega_\delta(\hat{S}_n) \le \varepsilon/2, T'_1 >n )
\P_x \Bigl(\max_{0 \le \ell \le b_n} \frac{|S_\ell|}{\sigma \sqrt{n}} \le \frac{\varepsilon}{4}, \tau_B > T'_{i-1},  T'_{i-1} \le b_n, H'_{i-1} \in dy \Bigr) \\
&\succeq & \frac{\sqrt{2} }{\sigma\sqrt{\pi n}} \inf_{r \le y \le a_n \sqrt{n}} \P_y \bigl( \omega_\delta(\hat{S}_n) \le \varepsilon/2 \, \bigl  | \, \bigr. T'_1 >n \bigr) \cdot \E_x (H_i' -H'_{i-1}) \I_{\{\kappa' \ge i, H'_{i-1} \in B_+ \}} .
\eeaa
The infimum vanishes since the family $\{ Law_y ( \hat{S}_n (\cdot)| T'_1 >n ) \}_{r \le y \le a_n \sqrt{n}, \, n \ge 1}$ is tight by \eqref{meander}, hence
summing over $i$ from $1$ to $k$ and letting $k \to \infty$, we conclude \eqref{tightness} by \eqref{main B} and \eqref{truncation}.

\underline{Distribution of the sign $\rho$.} We need to show that for $x \notin B$,
$$V_B(x) + x - \E_x S_{\tau_B} = 2\E_x \sum_{i =1}^{\kappa'} \bigl| H_{i}' - H_{i-1}' \bigr| \I_{\{H_{i-1}' \in B_+\}}.$$
If $B$ is an interval, this is immediate by $$V_B(x) = \E_x \Bigl [ \sum_{i =1}^{\kappa'}  (H_{i-1}' - H_{i}') \I_{\{H_{i-1}' \in B_+\}} +  \sum_{i =1}^{\kappa'} (H_{i}' - H_{i-1}' ) \I_{\{H_{i-1}' \in B_-\}} \Bigr ] =: \E_x [\Sigma_+ + \Sigma_-] $$ and
$\Sigma_+ - \Sigma_- = x - S_{T_{\kappa'}'} = x - S_{\tau_B}$. Neither of the two latter equations hold for a general $B$, in which case it should be proved that for any $x \notin B$,
$$\E_x \Bigl[ \sum_{i =1}^{\kappa'} ( S_{T_{i-1}'} - S_{T_i'}) \I_{\{S_{T_{i-1}'} \in Conv(B) \}} + (S_{T_{\kappa'}'} -S_{\tau_B}) \Bigr]= 0.$$

Define $\nu(0):=0$ and $\nu(k+1):= \min\{i > \nu(k): S_{T_i'} \in Conv(B)\}$ for $k \ge 0$. By conditioning on $S_{T_{\nu(1)}'}$ if $x \notin Conv(B)$, the required follows once we check that for any $y \in Conv(B) \setminus B$,
\be \label{stoppings}
\E_y \Bigl[ \sum_{k=0}^\infty S_{T_{\nu(k)}' \wedge \tau_B} - S_{T_{\nu(k)+1}' \wedge \tau_B} \Bigr] = 0.
\ee
Since all $T_{\nu(k)}' \wedge \tau_B$ and $T_{\nu(k)+1}' \wedge \tau_B$ are stopping times, the expectation of each term equals zero. Indeed, for the first term, $\E_y S_{T_1' \wedge \tau_B} = y$ by the optional stopping theorem, since $\E_y |S_{T_1' \wedge \tau_B} | \le \E_y (|S_{T_1'}| +  |S_{\tau_B} |) < \infty$ and
$$|\E_y S_n \I_{ \{ T_1' \wedge \tau_B >n \} } | \le (\E_y S_n^2 \cdot \P_y(T_1' \wedge \tau_B >n ))^{1/2} \le \sigma \sqrt{n} \P_y^{1/2}( T_1' >n ) \to 0$$ by \eqref{two sided}. Then, conditioning on $S_{T_{\nu(k)}' \wedge \tau_B}$ for $k \ge 1$, we see by induction that the expectations of the other terms are zero.

Finally, we have that
\beaa
\E_y \sum_{k=0}^\infty \bigl| S_{T_{\nu(k)}' \wedge \tau_B} - S_{T_{\nu(k)+1}' \wedge \tau_B} \bigr| &\le& (r-l + \sup_{l \le x \le r } \E_x |S_{T_1'}|), \E_y (\kappa' + 1) < \infty
\eeaa
where $\kappa'$ is integrable by the tail estimate $\P_y(\kappa' \ge i) \le \gamma'^{\lfloor i/3 \rfloor}$, see \eqref{Jump over'}. This allows to conclude the proof of \eqref{stoppings} by taking the expectation inside the sum.

The remaining statements of Theorem~\ref{THM MEANDER} for the case that $x=x_n \to \pm \infty$ are clear by $\P_{x_n}(\tau_B>n) \sim \P_{x_n}(T_1'>n)$.
\ep

\section{The largest gap: proofs} \label{SEC PROOFS GAP}

For any $\mathcal{R} \subset \R$, denote
\be \label{gap}
G(\mathcal{R}):= \sup \Bigl \{h \ge 0: \exists \, x \in \mathcal{R} \mbox{ such that } (x, x+h) \cap \mathcal{R} = \varnothing, x < \sup\{y: y \in \mathcal{R}\} \Bigr\}.
\ee
It is readily seen that $G(\mathcal{R})$ equals the largest spacing within the elements of at most countable $\mathcal{R}$ whenever its element admit ascending ordering, and in particular, $G_n=G(\{S_k\}_{k \ge 1}^n)$ and $G_- = G\bigl(\bigl\{S_n^{\geqslant}, -S_n^< \bigr\}_{n \ge 0} \bigr)$. Indeed, the left end of the largest spacing is an element of $\mathcal{R}$ while the condition $x < \sup\{y: y \in \mathcal{R}\}$ ensures that $(x, x+h)$ is within $\mathcal{R}$ (inside its convex hull). Thus \eqref{gap} can be taken as a definition of the largest gap, and it is the key formula to work with.

\bp[{\bf Proof of Proposition~\ref{THM GAP}}]
For any $h>0$, we have
\beaa
\{ G_n \ge 2h \} &= &\bigcup_{k =1}^n \bigl \{ S_i - S_k \notin (0,2h), 1 \le i \le n \bigr \}  \cap \bigl \{ S_k < M_n \bigr \} \\
&=& \bigcup_{k =1}^n \bigl \{ S_{k-i} - S_k \notin (0,2h), 1 \le i \le k-1 \bigr \}  \cap \bigl \{S_{k+i} - S_k \notin (0,2h), 1 \le i \le n-k \bigr \} \\
&& \cap \Bigl( \bigl \{ \max_{1 \le i \le k-1} (S_{k-i} - S_k) \ge 2h \bigr \} \cup  \bigl \{ \max_{1 \le i \le n-k} (S_{k+i} - S_k) \ge 2h \bigr \} \Bigr).
\eeaa
Reversing the time and changing the sign to account the trajectory of the walk by the time $k$, we obtain
\beaa
\P(G_n \ge 2h ) &\le& \sum_{k=1}^n \P_0(\tau_{(-2h, 0)} > k-1, T_1 \le k-1) \P_0(\tau_{(0,2h)} > n-k)\\
&+& \sum_{k=1}^n  \P_0(\tau_{(-2h, 0)} > k-1) \P_0(\tau_{(0,2h)} > n-k, T_1 \le n-k).
\eeaa
This implies that
\bea
\P(G_n \ge 2h ) &\le& \sum_{k=1}^n \bigl [ p_{k-1}^{(h)}(h) -  \P_h(T_1 > k-1 ) \bigr] p_{n-k}^{(h) }(-h) + p_{k-1}^{(h)}(h) \bigr[ p_{n-k}^{(h) }(-h) - \P_{-h}(T_1 > n-k)  \bigr] \notag \\
&=& \frac{2}{\sigma^2 \pi} \sum_{k=1}^n \frac{\bigl( V_h(h) - U_h(h)\bigr) V_h(-h) + \bigl(V_h(-h) - U_h(-h)\bigr) V_h(h) + o(1)}{\sqrt{k(n-k+1)}}, \label{beta}
\eea
and we get $$\limsup_{n \to \infty } \P(G_n \ge 2h ) \le \frac{2 }{\sigma^2 } \Bigl[ \bigl( V_h(h) - U_h(h)\bigr) V_{(0,2h)}(0) + \bigl(V_h(-h) - U_h(-h)\bigr) V_{(-2h,0)}(0) \Bigr].$$ Since $V_{(0,2h)}(0)$ and $V_{(-2h,0)}(0)$ are decreasing in $h$, the r.h.s. tends to zero as $h \to \infty$ by Theorem~\ref{THM EQUIV}, and the required follows.
\ep

Let us explain  the proof of Theorem~\ref{THM GAP LIMIT} given below. The key observation is \eqref{Int gap}, which means that gaps in the bulk of the range vanish. This fact is explained by the following heuristics. Suppose for some $1 \le k \le n$ and $h > u$, the interval $(S_k, S_k+h)$ is not hit by the walk by the time $n$ and $S_k$ is in the bulk of the range. Then $k$ splits the trajectory in two independent walks (one of them time-reversed) as used in the proof of Proposition~\ref{THM GAP}. Both trajectories avoid $(0,h)$ and by Theorem~\ref{THM MEANDER}, resemble Brownian meanders up to some possible fluctuations near $k$. Note that the meanders have different signs otherwise $S_k$ would be near the global extrema of the walk. Thus $k$ resembles a point of increase (or decrease) of a random walk. It is known that a centred random walk asymptotically has no points of increase, which contradicts our assumption on the existence of a gap in the bulk. The author owes to Yuval Peres this idea of referring to points of increase.

\bp[{\bf Proof of Theorem~\ref{THM GAP LIMIT}}]

\underline{No gaps in the bulk of the range.} Let us prove \eqref{Int gap}. For any $\varepsilon > u$,
\beaa
&&\bigl \{G_n^{Int} \ge \varepsilon \bigr \} \\
&=& \bigcup_{k=1}^n \biggl \{ \sum_{j=1}^n \I_{\{ S_j \le S_k\}} \ge b_n; \sum_{j=1}^n \I_{\{ S_j \ge S_k\}} \ge b_n; S_i \notin (S_k, S_k+ \varepsilon), 1 \le i \le n ; S_k < M_n \biggr \} \\
&\subset& \bigcup_{k=1}^n \biggl \{ \sum_{j=1}^{k-1} \I_{\{ S_{k-j} -S_k \le 0\}} \ge b_n/2 \, \bigvee \, \sum_{j=1}^{n-k} \I_{\{ S_{k+j} - S_k \le 0\}}  \ge b_n/2 ;  S_{k-i} -S_k \notin (0, \varepsilon), 1 \le i \le k-1; \\
&&\sum_{j=1}^{k-1} \I_{\{ S_{k-j} -S_k \ge 0\}} \ge b_n/2 \, \bigvee \,  \sum_{j=1}^{n-k} \I_{\{ S_{k+j} - S_k \ge 0\}}  \ge b_n/2;  S_{k+i} -S_k \notin (0, \varepsilon), 1 \le i \le n-k\biggr \}.
\eeaa
Notice that
$$\P \biggl ( \sum_{j=1}^{k-1} \I_{\{ S_{k-j} -S_k \le 0\}} \ge b_n/2 ; \sum_{j=1}^{k-1} \I_{\{ S_{k-j} -S_k \ge 0\}} \ge b_n/2 ;  S_{k-i} -S_k \notin (0, \varepsilon), 1 \le i \le k-1 \biggr ) = o(k^{-1/2})$$
holds uniformly in $k = k_n$ such that $b_n=o(k)$. This follows by Theorem~\ref{THM MEANDER} applied to the walk $-S_n$ and the fact that
$$\P \biggl ( \sum_{j=1}^{k-1} \I_{\{ S_j = 0\}} \ge b_n/4 ; -S_i \notin (0, \varepsilon), 1 \le i \le k-1 \biggr ) = o(k^{-1/2}).$$ Then, using the analogous observation that the probability that the other two sums simultaneously exceed $b_n/2$ equals $o((n-k)^{-1/2})$ for $k$ such that $b_n=o(n-k)$, we argue as in \eqref{beta} to obtain that
\beaa
\P(G_n^{Int} \ge \varepsilon ) &\le&
\P \biggl (\bigcup_{k=1}^n \biggl \{ \sum_{j=1}^{k-1} \I_{\{ S_{k-j} -S_k \ge 0\}} + \sum_{j=1}^{n-k} \I_{\{ S_{k+j} - S_k \le 0\}}  < b_n \biggr \} \biggl ) \\
&+& \P \biggl (\bigcup_{k=1}^n \biggl \{ \sum_{j=1}^{k-1} \I_{\{ S_{k-j} -S_k \le 0\}} + \sum_{j=1}^{n-k} \I_{\{ S_{k+j} - S_k \ge 0\}}  < b_n \biggr \} \biggl ) + o(1).
\eeaa

The second probability in the r.h.s. equals the first one computed for the walk $-S_n$, and thus \eqref{Int gap} follows if
we prove that
\be \label{no increase}
\lim_{n \to \infty} \P \biggl ( \min_{1 \le k \le n }  \biggl [ \sum_{j=1}^k \I_{\{ S_j > S_k \}} + \sum_{j=k}^n \I_{\{ S_j < S_k \}}  \biggr] < b_n \biggl ) = 0.
\ee

Let $f(t)$ be a continuous function on $[0,1]$. Define $$H[f](t):= \int_0^t \I_{\{ f(s) > f(t)\}} ds +  \int_t^1 \I_{\{ f(s) < f(t)\}} ds, \quad F[f]:= \inf_{0 \le t \le 1 } H[f](t).$$ The functions $\I_{\{ a > f(t)\}}$ and $\I_{\{ a < f(t)\}}$ are lower-semicontinuous for any fixed $a$. Then the Fatou lemma implies that $H[f]$ is lower-semicontinuous; moreover, by the dominated convergence theorem, the function $H[f]$ is continuous at any point of the set $$ C_f:=\{s \in [0,1]: \mbox{mes}(u \in [0,1]: f(u)=f(s)) = 0 \}.$$ Since $H[f]$ is lower-semicontinuous and $[0,1]$ is a compact set, the value $F[f]$ is attained at some
$t^* \in [0,1]$. In particular, this implies that $F[f]=0$ iff the function $f$ has a point of increase, i.e. a point $t$ such that $f(s) \le f(t)$ for any $s \le t$ and $f(s) \ge f(t)$ for any $s \ge t$.

We claim that $F[f]:C[0,1] \to [0,1]$ is continuous at any $f \in C[0,1]$ such that $C_f=[0,1]$. Indeed, let $f_n \to f$ in $C[0,1]$ and $t_n \to t$. Then $f_n(t_n) \to f(t)$, and $H[f_n](t_n) \to H[f](t)$ by the dominated convergence theorem since $\I_{\{ f_n(s) \gtrless f_n(t_n)\}} \to \I_{\{ f(s) \gtrless f(t)\}}$ a.s. Now, if $t_n^*$ are such that $F[f_n]= H[f_n](t_n^*)$, then for any converging subsequence $t_{n_k}^*$, it is true that $$\lim_{k \to \infty} F[f_{n_k}]= \lim_{k \to \infty} H[f_{n_k}](t_{n_k}^*) = H[f](\lim_{k \to \infty} t_{n_k}^*) \ge F[f],$$ hence $\liminf_n F[f_n] \ge F[f]$. On the other hand, $F[f] = H[f](t^*)$ for some $t^*$, and $$\limsup_{n \to \infty} F[f_n] \le \lim_{n \to \infty} H[f_n](t^*) = H[f](t^*) = F[f],$$ and the claim follows. {\it Nota bene:} the assumption $C_f=[0,1]$ is actually required only to prove the last equation.

Let $W$ be a standard Brownian motion. Now since $\hat{S}_n \stackrel{\D}{\to} W$ in $C[0,1]$ and $\P(C_W = [0,1]) = 1$, it holds that $F[\hat{S}_n] \stackrel{\D}{\to} F[W]$. By the theorem of Dvoretzky, Erd{\"o}s and Kakutani, $W$ has no points of increase (see, e.g. M{\"o}rters and Peres~\cite[Sec. 5.2]{MortersPeres}) so $F[W]$ has no atom at zero. Then \eqref{no increase} follows since $b_n /n \to 0$, and \eqref{Int gap} is proved.

\underline{Gaps near the minimum of the range.} For any $b \ge 1$, denote
$$I_\geqslant(b):=\Bigl \{k  \ge 0:  \sum_{i=0}^\infty \I_{\{ S_i^{\geqslant} \le S_k^{\geqslant}  \}} + \I_{\{ -S_i^< \le S_k^{\geqslant} \}} > b \Bigr \},$$
and introduce $I_<(b)$ by analogy, replacing $S_k^{\geqslant}$ by $-S_k^<$ in the definition of $I_\geqslant(b)$. Recall that by assumption, the Markov chains $S_k^{\geqslant}$ and $S_k^<$ are independent and start at $0$. Let
$$G(b):= G \bigl(\{ S_k^{\geqslant} \}_{k \notin I_\geqslant(b)} \cup \{-S_k^< \}_{k \notin I_<(b)} \bigr), \qquad G_n^-(b):= \max_{1 \le k \le b-1} S_{(k+1,n)} -S_{(k,n)}$$ be the largest gaps within the $b$ smallest values among $\{S_k^{\geqslant}, -S_k^<\}_{n \ge 0}$ and $\{S_k\}_{k=1}^n$, respectively. Notice that $G(b)$ is a proper random variable by $G(b) \le G \bigl(\{ S_k^{\geqslant} \}_{k < b} \bigr)$.

Denote by $\mu_n$ the location of the first global minimum of the walk by the time $n$. It follows from  (\ref{half-tail}$'$) and \eqref{Doob h-limit} that for any $\ell \ge 1$ and any bounded measurable $f: \R^{2 \ell + 1} \to \R$,
$$ \lim_{n \to \infty} \E \bigl( f \bigl( S_{\mu_n - \ell} -  S_{\mu_n}, \dots, S_{\mu_n+\ell}- S_{\mu_n} \bigr)  \bigl | \bigr. \mu_n = k_n \bigr) = \E_ 0 f \bigl(-S_\ell^<, \dots, -S_1^<, 0, S_1^{\geqslant}, \dots,   S_\ell^{\geqslant} \bigr)$$ uniformly in $k_n$ such that $k_n \wedge n-k_n \to \infty$; we will need only the unconditioned version of this equation, which is  weaker. On the other hand, for any fixed $b$,
$$\lim_{\ell \to \infty} \P_0 \bigr( I^c_\geqslant(b) \cup I^c_<(b) \subset [0, \ell] \bigl) = 1, \qquad \lim_{\ell \to \infty } \liminf_{n \to \infty}\P_0 \Bigl( \{S_{(i,n)}\}_{i \le b} \subset \{ S_{\mu_n + i} \}_{|i| \le \ell} \Bigr) = 1,$$ where the first equation follows by transience of the Markov chains  $S_n^{\geqslant}$ and $-S_n^<$, and the second holds by a result of Ritter~\cite{Ritter}: for any $\beta \in (0, 1/2)$,
$$\lim_{\delta \to 0} \liminf_{n \to \infty} \P_0 \bigl(S_k \ge \delta k^{\beta}, \, 1 \le k \le n \bigl| \bigr. \tau_{(-\infty, 0]} > n \bigr) = 1.$$
Recalling that $G(b)$ is a proper random variable, we obtain that for any $b \ge 1$,
\be \label{partial gaps}
G_n^-(b) \stackrel{\D}{\longrightarrow} G(b).
\ee

We claim that
\be \label{gaps =}
\lim_{b \to \infty} \liminf_{n \to \infty} \P \bigl(G_n^- (b) = G_n^-(b_n)\bigr) = 1.
\ee
Assume the converse. Then there exist sequences $b_i', n_i \to \infty$ such that
\be \label{absurdum}
\lim_{i \to \infty} \P \bigl(  G_{n_i}^- (b_i') \ge \max_{b_i' \le k \le b_{n_i}} S_{(k,n_i)} - S_{(k-1,n_i)} \bigr) <1 .
\ee
If the walk is $\lambda$-arithmetic, then \eqref{Int gap}, considered for the sequence $b_i'$, implies
$$\lim_{i \to \infty} \P \bigl(  \max_{b_i' \le k \le n_i - b_{n_i}'} S_{(k,n_i)} - S_{(k-1,n_i)} = \lambda \bigr) = 1,$$
which contradicts \eqref{absurdum} by $\P( G_{n_i}^- (b_i') \ge \lambda) = 1 -\P^{n_i}(X_1=0) \to 1$. If the walk is non-arithmetic, it follows from \eqref{partial gaps} that  for any $b, \varepsilon >0$,
$$\liminf_{i \to \infty} \P( G_{n_i}^- (b_i')  > \varepsilon) \ge \liminf_{i \to \infty} \P( G_{n_i}^- (b)  > \varepsilon) \ge \P(G(b) > \varepsilon).$$
Combined with \eqref{Int gap}, this contradicts \eqref{absurdum}. Indeed, the later probability can be made arbitrarily close to $1$ by choosing $\varepsilon$ to be small enough and then choosing a large enough $b$. Here we used that $\P(G > 0) = 1$, which holds by transience of $S_n^{\geqslant}$ and $-S_n^<$, and $G(b) \nearrow G$ a.s.

The family $\{ G_n^-(b_n) \}_{n \ge 1}$ is tight by the fact that $G_n^-(b_n) \le G_n$ and Proposition~\ref{THM GAP}. Then there exists a weakly converging subsequence, and its weak limit is exactly $G$ by \eqref{partial gaps}, \eqref{gaps =} and the reverse of Theorem 3.2 by Billingsley~\cite{Billingsley}. Then $G$ is a proper random variable, hence $G(b) \stackrel{\D}{\longrightarrow} G$, and thus $G_n^-(b_n) \stackrel{\D}{\longrightarrow} G$ by the direct theorem.

\underline{The limit distribution.} An analogous consideration of gaps near the last global maximum $\eta_n$ implies that $G_n^+(b_n) \stackrel{\D}{\longrightarrow} G$, where $G_n^+(b):= \max_{1 \le k \le b-1} S_{(k+1,n)} -S_{(k,n)}$. Now \eqref{Ext gap} follows by $G_n^{Ext} = \max(G_n^-(b_n), G_n^+(b_n))$ and \eqref{gaps =} if we show that $G_n^-(b)$ and $G_n^+(b)$ are asymptotically independent for any fixed $b \ge 1$.

Arguing analogously to the proof of \eqref{partial gaps}, it suffices to check that for any $\ell \ge 1$ and any bounded measurable $f, g: \R^{2 \ell + 1} \to \R$,
\bea
&&\lim_{n \to \infty} \E \bigl( f \bigl( S_{\mu_n - \ell} -  S_{\mu_n}, \dots, S_{\mu_n+\ell}- S_{\mu_n} \bigr) g \bigl( S_{\eta_n} - S_{\eta_n - \ell}, \dots, S_{\eta_n} - S_{\eta_n+\ell} \bigr) \bigl | \bigr. \mu_n = k_n, \eta_n = m_n \bigr) \notag \\
&& = \E_ 0 f \bigl(-S_\ell^<, \dots, -S_1^<, 0, S_1^{\geqslant}, \dots,   S_\ell^{\geqslant} \bigr) \E_0 g \bigl(S_\ell^{\geqslant}, \dots,   S_1^{\geqslant}, 0, -S_1<, \dots, -S_\ell^< \bigr) \label{asymp indep}
\eea
uniformly in $k_n$ and $m_n$ such that $k_n \wedge m_n \wedge  n -  k_n \wedge n - m_n   \wedge |m_n - k_n| \to \infty$. W.l.o.g., assume that $k_n < m_n$.

For any $1 \le k \le n$, denote
$$E_k^-:= \frac{1}{\sqrt{k}} \max_{1 \le i \le k} S_i -  S_k,  \quad E_k^+:=  \frac{1}{\sqrt{n-k}} \max_{ k  \le i \le n} S_k - S_i$$
(where $E_n^+:=1$) and for any $1 \le k < m \le n$, denote
$$I_{k , m}^-:= \frac{1}{\sqrt{(m - k)/2}} \max_{k \le i \le \lfloor (m+k)/2 \rfloor} S_i -  S_k, \quad I_{k , m}^+:= \frac{1}{\sqrt{(m - k)/2}} \max_{ \lfloor (m+k)/2 \rfloor \le i \le m} S_m -  S_i$$ and
$$D_{k , m}^-:= \frac{S_{\lfloor (m+k)/2 \rfloor} -  S_k}{\sqrt{(m - k)/2}}, \quad D_{k , m}^+:= \frac{S_m -  S_{\lfloor (m+k)/2 \rfloor }}{\sqrt{(m - k)/2}}.$$

We first study the probability of the condition in \eqref{asymp indep}. By combining \eqref{half-tail}, $(\ref{half-tail}')$ and \eqref{meander} with the fact that $\E_0 H_+' \cdot \E_0 H_- = \sigma^2/2$ (see Feller~\cite[Sec. XVIII.4]{Feller}), we obtain
\bea
&& \P(\mu_n = k_n, \eta_n = m_n) \times \bigl( \pi (m_n - k_n) \sqrt{k_n(n-m_n)} \bigr)  \notag \\
&\sim& \P \biggl( \max \Bigl( \sqrt{\frac{2k_n}{m_n-k_n}}E_{k_n}^-,  I_{k_n , m_n}^-,  I_{k_n , m_n}^+, \sqrt{\frac{2(n-m_n)}{m_n-k_n}}E_{m_n}^+ \Bigr) \le D_{k_n , m_n}^- + D_{k_n , m_n}^+ \biggr) \notag \\
&\sim& \P \biggl( \max \Bigl( \sqrt{\frac{2k_n}{m_n-k_n}} \overline{W}_+^{(1)},  \overline{W}_+^{(2)},  \overline{W}_+^{(3)}, \sqrt{\frac{2(n-m_n)}{m_n-k_n}} \overline{W}_+^{(4)} \Bigr) \le W_+^{(2)}(1) + W_+^{(3)}(1) \biggr) \label{condit}.
\eea
where $W_+^{(i)}(\cdot)$ are independent standard Brownian meanders on $[0,1]$ and $\overline{W}_+^{(i)}:= \max_{0 \le t \le 1} W_+^{(i)}$. Let us proceed analogously with the l.h.s. of \eqref{asymp indep}. We start conditioning on $X_{k_n - \ell +1}, \dots, X_{k_n + \ell},$ $X_{m_n - \ell +1}, \dots, X_{m_n + \ell}$. We assume that $n$ is large enough such that $m_n -k_n > 2 \ell$. Consider, for example, the part of trajectory between $m_n$ and $n$, which forms a new random walk $\tilde{S}_i:= S_{m_n+i} - S_{m_n}$. The first $\ell$ steps are already negative by fixing appropriate values of $X_{m_n +1}, \dots, X_{m_n + \ell}$. The key observation is that by \eqref{meander}, the remaining part of the trajectory of $\tilde{S}_i$ with $\ell + 1 \le i \le n-m_n$ conditioned on $\{\tau_{[0, \infty)}' > n -m_n - \ell \}$ under $P_{\tilde{S}_\ell}$, converges to the same meander $W_+^{(4)}(\cdot)$ as the whole trajectory of $\tilde{S}_i$ with $1 \le i \le n-m_n$ conditioned on $\{ \tau_{[0, \infty)}' > n -m_n \}$ under $P_0$. We conclude the proof of \eqref{asymp indep} by making analogous observations for the other three parts of the trajectory of $S_n$ and recalling \eqref{condit}.


\underline{Convergence of the number of non-visited sites.} Here we prove \eqref{empty}.

W.l.o.g., assume  that the random walk $S_n$ is $1$-arithmetic. By \eqref{Int gap} and \eqref{Ext gap} it holds that $\P(E_n \le 2 b_n (G_n -1)) \to 1$ for any $b_n \to \infty$, hence $E_n/b_n' \stackrel{\P}{\to} 0$ for any $b_n' \to \infty$ by the tightness of $Law(G_n)_{n \ge 1}$. This implies that $Law(E_n)_{n \ge 1}$ is tight. Indeed, assume the converse. Then there exist increasing sequences $n_k, r_k \to \infty$ such that $\liminf_k \P(E_{n_k} > r_k) >0$, which contradicts $E_n/b_n' \stackrel{\P}{\to} 0$ with the sequence $b_n'$ defined by $b_n':=r_k$ for $n_k \le n < n_{k+1}$.

It remains to argue analogously to the proof of \eqref{Ext gap} to conclude \eqref{empty}.
\ep

\section*{Acknowledgements}
The author thanks Vitali Wachtel and Yuval Peres for sharing ideas and useful discussions, and Ron Doney for comments and his interest to this paper.


\begin{thebibliography}{99}


\bibitem{Baumgarten} Baumgarten C. (2013+) Persistence of some iterated processes, preprint.

\bibitem{Belkin} Belkin B. (1972) An invariance principle for conditioned recurrent random walk attracted to a stable law. {\it Z. Wahrsch. Verw. Gebiete} \textbf{21} 45-–64.

\bibitem{BertoinDoney} \textsc{Bertoin, J.} and \textsc{Doney, R. A.} (1994) On conditioning a random walk to stay nonnegative. \textit{Ann. Probab.} \textbf{22} 2152--2167.

\bibitem{Billingsley} \textsc{Billingsley, P.} (1999) \textit{Convergence of Probability Measures, 2nd Edition}. Wiley.

\bibitem{Borodin} Borodin A.N. (1982) On the asymptotic behavior of local times of recurrent random walks with finite variance. {\it Theory Probab. Appl.} \textbf{26} 758–-772.

\bibitem{BGT} \textsc{Bingham, N. H.}, \textsc{Goldie, C. M.} and \textsc{Teugels, J. L.} (1989) \textit{Regular Variation.} Cambridge University Press.

\bibitem{Doney} \textsc{Doney, R. A.} (2010) Local behaviour of first passage probabilities. \textit{Probab. Theory Relat. Fields} \textbf{152} 559--588.

\bibitem{Eppel} \textsc{Eppel, M. S.} (1979) A local limit theorem for the first overshoot. \textit{Siberian Math. J.} \textbf{20} 130--138.

\bibitem{Gut} \textsc{Gut, A.} (2009) \textit{Stopped Random Walks: Limit Theorems and Applications, 2nd ed.} Springer.

\bibitem{Feller} \textsc{Feller, W.} (1966) \textit{An Introduction to Probability Theory and Its Applications, Vol. 2}. Wiley, New York.

\bibitem{KestenSpitzer} \textsc{Kesten, H.} and \textsc{Spitzer, F.} (1963) Ratio theorems for random walks I. \textit{J. Anal. Math.} \textbf{11} 285--322.

\bibitem{MortersPeres} M{\"o}rters P. and Peres Y. (2010) \textit{Brownian motion.} Cambridge University Press.

\bibitem{Nagaev} \textsc{Nagaev, A. V.} (1976) Some limit theorem of the renewal theory. \textit{Theor. Probab. Appl.} \textbf{20} 323--336.

\bibitem{Problems} Polya G. and Szego G. (1997) {\it Problems and Theorems in Analysis I, reprint.} Springer.

\bibitem{Ritter} Ritter, G. A. (1981) Growth of random walks conditioned to stay positive. \textit{Ann. Probab.} \textbf{9} 699--704.

\bibitem{Spitzer} \textsc{Spitzer, F.} (2001) \textit{Principles of Random Walk, 2nd ed.} Springer.

\bibitem{MeIterated} Vysotsky V. (2013+) Persistence of iterated random walks, in progress.


\end{thebibliography}
\end{document}